\documentstyle[11pt,amstex,amssymb]{amsart}
\newtheorem{thm}{Theorem}[section]
\newtheorem{prop}[thm]{Proposition}
\newtheorem{lemma}[thm]{Lemma}

\renewcommand{\thethm}{\arabic{section}.\arabic{thm}}
\renewcommand{\theequation}{\thesection.\arabic{equation}}
\newcommand{\proof}{\noindent{\it Proof.}\enspace}

\def\bR{{\bold R}}
\def\bT{{\bold T}}

\def\cX{{\cal X}}
\def\cM{{\cal M}}
\def\Ric{{\rm Ric}}
\def\Hess{{\rm Hess}}
\def\eps{\varepsilon}
\def\bN{{\bold N}}
\def\Re{{\rm Re}\,}

\def\Tr{{\rm Tr}}
\def\U{{\rm U}}
\def\SU{{\rm SU}}
\def\SO{{\rm SO}}
\def\i{\sqrt{-1}}

\begin{document}

\title[Inequalities related to free entropy]
{Inequalities related to free entropy derived from 
random matrix approximation}
\author[F. Hiai]{Fumio Hiai$\,^1$}
\address{Graduate School of Information Sciences,
Tohoku University, Aoba-ku, Sendai 980-8579, Japan}
\author[D. Petz]{D\'enes Petz$\,^2$}
\address{Department for Mathematical Analysis,
Budapest University of Technology and Economics,
H-1521 Budapest XI., Hungary}
\author[Y. Ueda]{Yoshimichi Ueda$\,^3$}
\address{Graduate School of Mathematics,
Kyushu University, Fukuoka 810-8560, Japan}

\thanks{$^1\,$Supported in part by
Grant-in-Aid for Scientific Research (C)14540198 and by Strategic
Information and Communications R\&D Promotion Scheme of MPHPT}

\thanks{$^2\,$Supported in part by MTA-JSPS project (Quantum Probability
and Information Theory) and by OTKA T032662.}

\thanks{$^3\,$Supported in part by Grant-in-Aid for Young
Scientists (B)14740118.}

\maketitle

\begin{abstract}
Biane proved the free analog of the logarithmic Sobolev inequality for
probability measures on $\bR$ by means of random matrix approximation
procedure. We show that the same method can be applied to reprove Biane and
Voiculescu's free analog of Talagrand's transportation cost inequality for
measures on $\bR$. Furthermore, we prove the free analogs of the logarithmic
Sobolev inequality and the transportation cost inequality for measures on
$\bT$ as well by extending the method to special unitary random matrices. 
\end{abstract}

\section*{Introduction}

Since  its first systematic study done by L.\ Gross \cite{G} in
1975, the logarithmic Sobolev inequality (LSI) has been discussed by many
authors in various contexts, in particular, in close connection with the
notions of hypercontractivity and spectral gap. An LSI can be
understood to compare the relative Fisher information with the relative
entropy. Among other things, we here refer to the LSI due to D.~Bakry and
M.~Emery \cite{BE} in the general Riemannian manifold setting, which is of
quite use for our present purpose. Another interesting inequality was
presented by M.~Talagrand \cite{Ta} in 1996, called the transportation cost
inequality (TCI). A TCI compares the (quadratic) Wasserstein distance
$W(\mu,\nu)$ between probability measures $\mu,\nu$ (for the definition see
\eqref{F-4.1} in \S4 of this paper) with $\sqrt{S(\mu,\nu)}$, the square
root of the relative entropy. Indeed, in \cite{Ta} Talagrand proved the
inequality $W(\mu,\nu)\le\sqrt{S(\mu,\nu)}$ when $\nu$ is the standard
Gaussian measure on $\bR^n$, and an exposition in the case of more general
$\nu$ can be found in \cite{Le} for example. On the other hand, in \cite{OV}
F.~Otto and C.~Villani succeeded in discovering links between the LSI and the
TCI in the Riemannian manifold setting. This, combined with \cite{BE},
implies the TCI in the same situation as Bakry and Emery's LSI. See
\cite{Le2, Le, Vi} for more about these classical LSI and TCI as well as
related topics.

The relative free entropy $\widetilde\Sigma_Q(\mu)$ and the relative free
Fisher information $\Phi_Q(\mu)$ were introduced by Ph.\ Biane and
R.\ Speicher \cite{BS} for $\mu\in\cM(\bR)$, the probability measures on
$\bR$, relative to a real continuous function $Q$ on $\bR$, where $Q$ has a
certain growth in the case of $\widetilde\Sigma_Q(\mu)$ and it is a $C^1$
function in the case of $\Phi_Q(\mu)$. Note that $\widetilde\Sigma_Q(\mu)$
is regarded as the relative version of the free entropy $\Sigma(\mu)$
introduced by D.~Voiculescu \cite{V1} as the classical relative entropy is
the relative version of the Boltzmann-Gibbs entropy, while $\Phi_Q(\mu)$ in
the case $Q\equiv0$ reduces to the free Fisher information $\Phi(\mu)$ in
\cite{V1}. (The ``free relative entropy" $\Sigma(\mu,\nu)$ for two measures
was introduced in \cite{HMP} from a slightly different viewpoint.) In this
paper we introduce the relative free entropy $\widetilde\Sigma_Q(\mu)$ and
the relative free Fisher information $F_Q(\mu)$ for $\mu\in\cM(\bT)$, the
probability measures on the unit circle $\bT$, as well relative to a
real continuous function $Q$ on $\bT$ (being a $C^1$ function for
$F_Q(\mu)$). When $Q\equiv0$ the quantity $F_Q(\mu)$ becomes the free Fisher
information $F(\mu)$ introduced by Voiculescu \cite{V6}. An important fact
is that the relative free entropy $\widetilde\Sigma_Q(\mu)$ is the rate
function (or the so-called weighted logarithmic integral up to an additive
constant) of a large deviation for the empirical eigenvalue distribution of
a certain random matrix. Indeed, $\widetilde\Sigma_Q(\mu)$ for
$\mu\in\cM(\bR)$ is the good rate function of large deviation principle for
the $n\times n$ selfadjoint random matrix determined by the function $Q$,
while $\widetilde\Sigma_Q(\mu)$ for $\mu\in\cM(\bT)$ is that for the
$n\times n$ (special) unitary random matrix associated with $Q$. The
definitions of these quantities as well as related matters are collected in
the first \S1 of this paper.

Voiculescu's inequality in \cite[Proposition 7.9]{V5} is the first free
probabilistic analog of the LSI. Extending its single variable case (see
\eqref{Voiculescu's LSI} in \S2), Biane obtained in \cite{Bi} the following
free LSI:
$$
\widetilde\Sigma_Q(\mu)\le{1\over2\rho}\Phi_Q(\mu)
\quad\text{for $\mu\in\cM(\bR)$}
$$
if $Q''(x)\ge\rho$ on $\bR$ with a constant $\rho>0$. To prove this, Biane
applied the classical LSI on the Euclidean space to the related selfadjoint
random matrices as mentioned above and used the weak convergence of their
mean eigenvalue distributions. Although the differentiability assumption of
$Q$ is not quite explicitly written in \cite{Bi}, Biane's free LSI is
certainly valid if $Q$ is a $C^1$ function such that $Q(x)-{\rho\over2}x^2$
is convex on $\bR$. For the sake of completeness, in \S2 we give a proof of
this general case by a usual approximation technique.

The first main aim of this paper is to show the variant of Biane's free LSI
for measures on $\bT$. In \S3 we prove
$$
\widetilde\Sigma_Q(\mu)\le{1\over1+2\rho}F_Q(\mu)
\quad\text{for $\mu\in\cM(\bT)$}
$$
if $Q$ is a $C^1$ function on $\bT$ such that
$Q\Bigl(e^{\i t}\Bigr)-{\rho\over2}t^2$ is convex on $\bR$ with a constant
$\rho>-1/2$. The proof is based on random matrix approximation. We can
apply Bakry and Emery's classical LSI on the special unitary group $\SU(n)$,
a Riemannian manifold, to the related $n\times n$ special unitary random
matrices and pass to the scaling limit as $n$ goes to $\infty$. Here, we
need the convergence of the empirical eigenvalue distribution of the random
matrix not only in the mean but also in the almost sure sense that is a
consequence of the corresponding large deviation principle. Although the
large deviation theorem (Theorem \ref{T-1.2} below) for ``special" unitary
random matrices is essentially same as that for unitary random matrices
shown in \cite{HP2}, the proof is a bit more complicated so that we sketch
it in Appendix for the convenience of the reader. We also need a few stuffs
from differential geometry, in particular, the exact computation of the
Ricci curvature tensor of $\SU(n)$ (with respect to the Riemannian structure
associated with the usual trace on $M_n({\mathbf C})$) to check the so-called
Bakry and Emery criterion (see \S\S1.7).

In \cite{BV} Biane and Voiculescu obtained the free analog of Talagrand's
TCI for compactly supported $\mu\in\cM(\bR)$ as follows:
$$
W(\mu,\gamma_{0,2})
\le\sqrt{-\Sigma(\mu)+\int{x^2\over2}\,d\mu(x)-{3\over4}},
$$
where $\gamma_{0,2}$ denotes the standard semicircular distribution (with
radius $2$). Their proof involves the free process and the complex Burgers'
equation, and it is a realization of free probability parallel of not only
the result itself but also the proof in \cite{OV}. The proof itself
justifies the above inequality to be the right free analog of Talagrand's
TCI.

Our second main aim is to reprove Biane and Voiculescu's free TCI in a
slightly more general setting by making use of random matrix approximation
and furthermore to give a free TCI for measures on $\bT$ in a similar way.
This aim is our initial motivation; we first wanted to find another proof to
Biane and Voiculescu's TCI by use of random matrix approximation on the
lines of so-called Voiculescu's heuristics in \cite{V1} and to justify
Biane and Voiculescu's TCI as the right free analog from the viewpoint of
random matrix theory. In \S4 we prove the free TCI
$$
W(\mu,\mu_Q)\le\sqrt{{1\over\rho}\widetilde\Sigma_Q(\mu)}
\quad\text{for compactly supported $\mu\in\cM(\bR)$}
$$
if $Q$ is a real function on $\bR$ such that $Q(x)-{\rho\over2}x^2$ is
convex with a constant $\rho>0$ and $\mu_Q$ is the equilibrium measure
associated with $Q$ (or the unique minimizer of $\widetilde\Sigma_Q(\mu)$
for $\mu\in\cM(\bR$)). When $Q(x)=x^2/2$ and $\rho=1$, this becomes Biane
and Voiculescu's TCI. To prove this, we first suppose that $\mu$ is
supported in $[-R,R]$ and that $Q_\mu(x):=2\int_\bR\log|x-y|\,d\mu(y)$ is
continuous on $\bR$. We consider two $n\times n$ selfadjoint random
matrices; one is associated with $Q$, and the other is associated with
$Q_\mu$ and restricted on the $n\times n$ selfadjoint matrices with the
operator norm $\le R$. Then, these random matrices are probability measures
on the space of $n\times n$ selfadjoint matrices ($\cong\bR^{n^2}$), and
the classical TCI for these measures asymptotically approaches, as $n$
goes to $\infty$, to the free TCI we want. The case of general compactly
supported $\mu\in\cM(\bR)$ can be treated by an approximation technique.
Furthermore, as presented in \S5, a similar method using special unitary
random matrices can work to prove the free TCI
$$
W(\mu,\mu_Q)\le\sqrt{{2\over1+2\rho}\widetilde\Sigma_Q(\mu)}
\quad\text{for $\mu\in\cM(\bT)$}
$$
if $Q$ is such a real function on $\bT$ as in the free LSI, that is,
$Q\Bigl(e^{\i t}\Bigr)-{\rho\over2}t^2$ is convex on $\bR$ with $\rho>-1/2$.
Here, $W(\mu,\mu_Q)$ is the Wasserstein distance with respect to the
geodesic distance (or the angular distance) on $\bT$. In the particular
case where $Q\equiv0$ and $\rho=0$, we have
$W(\mu,d\theta/2\pi)\le\sqrt{2\Sigma(\mu)}$.

In this way, we clarify the advantage of random matrix approximation
procedure in studying free probabilistic analogs of certain classical
theories involving relative entropy and/or Fisher information. The present
paper may be regarded as one more attempt subsequent to \cite{BG, HMP}
toward rigorous realizations of Voiculescu's heuristics in
\cite{V1} which claims that the classical entropy of random matrices, if
suitably arranged, asymptotically converges to the free entropy of the
limit distribution as the matrix size goes to infinity.

The final \S6 is a collection of remarks, examples and related results; in
particular, we give the variants of the above free LSI and TCI for measures
on the half line $\bR^+$.

\section{Preliminaries}
\setcounter{equation}{0}

The purpose of this preliminary section is to summarize, for the convenience
of the reader, the basic notions and the results which will be needed
later. We will use them with no explicit explanation in the main part of
this paper.   

\subsection{Notations} The set of all Borel probability measures on a
Polish space ${\mathcal X}$ is denoted by ${\mathcal M}({\mathcal X})$. 
The Dirac measure at a point $x\in{\mathcal X}$ is denoted by $\delta_x$ as
usual. For $\mu, \nu \in {\mathcal M}({\mathcal X})$, the {\it relative
entropy} of $\mu$ with respect to $\nu$ is denoted by $S(\mu,\nu)$, which is
defined by 
\begin{equation}\label{relative entropy}
S(\mu,\nu) := \int_{\mathcal X} \log\frac{d\mu}{d\nu}\,d\mu
=\int_{\mathcal X} \frac{d\mu}{d\nu}\log\frac{d\mu}{d\nu}\,d\nu 
\end{equation}
when $\mu$ is absolutely continuous with respect to $\nu$; otherwise
$S(\mu,\nu) := +\infty$.  

The usual trace on $M_n({\mathbf C})$, the $n\times n$ complex matrices, is
denoted by $\Tr_n$. The Hilbert-Schmidt norm on $M_n({\mathbf C})$ induced
from ${\mathrm Tr}_n$ is denoted by $\Vert \cdot \Vert_{HS}$, i.e.,
$\Vert A\Vert_{HS} := {\mathrm Tr}_n(A^* A)^{1/2}$ for
$A \in M_n({\mathbf C})$.  Let $M_n^{sa}$ denote the set of all $n\times n$
self-adjoint matrices, ${\mathrm U}(n)$ the group of all $n\times n$
unitaries, and ${\mathrm SU}(n)$ the special unitary group of order $n$,
i.e., the group of all $n\times n$ unitaries whose determinants are equal to
one.

\subsection{Free entropy and free Fisher information for measures on
${\mathbf R}$}\label{Free entropy on R}
The notions of free entropy and free Fisher information are the free
probabilistic analogs of the Boltzmann-Gibbs entropy and the Fisher
information in classical theory. For each $\mu\in\cM(\bR)$, Voiculescu
\cite{V1} introduced the {\it free entropy} of $\mu$
$$
\Sigma(\mu) := \iint_{{\mathbf R}^2} \log\left|x-y\right|
d\mu(x)\,d\mu(y), 
$$
which is the minus of the so-called {\it logarithmic energy} of $\mu$
useful in potential theory (see \cite{ST}). It is the ``main component" of
the free entropy $\chi(\mu)$ of $\mu$ introduced in \cite{V2}:
\begin{equation}\label{free entropy chi}
\chi(\mu) = \Sigma(\mu) + \frac{3}{4} + \frac{1}{2}\log 2\pi.
\end{equation}

Assume that $\mu \in {\mathcal M}({\mathbf R})$ has the density
$p = d\mu/dx$ (with respect to the Lebesgue measure $dx$) belonging
to the $L^3$-space $L^3({\mathbf R}) := L^3\left({\mathbf R},dx\right)$.
In \cite{V1} Voiculescu also introduced the {\it free Fisher information}
of $\mu$  
$$
\Phi(\mu) := \frac{4\pi^2}{3} \int_{\mathbf R} p(x)^3\,dx =
\frac{4\pi^2}{3} \Vert p \Vert_3^3.  
$$
The {\it Hilbert transform} of $p$   
\begin{equation}\label{Hilbert transf on reals}
(Hp)(x) := \lim_{\varepsilon \searrow 0} 
\int_{|x-t|>\varepsilon} \frac{p(t)}{x-t}\,dt 
\end{equation} plays an important role in the study of free Fisher
information. The limit in \eqref{Hilbert transf on reals} really exists for
a.e.\ $x \in {\mathbf R}$ (as long as $p\in L^q(\bR)$ with $1 < q < \infty$),
and it is known that $p \in L^q(\bR)$ implies $Hp \in L^q(\bR)$ for each
$1 < q < \infty$. See \cite[Chapter VI]{Ko} for the Hilbert transform on
$\bR$. As shown in \cite[Lemma 3.3]{V1} we see that 
\begin{equation}\label{integral formula for Hilbert transf}
\int_{\mathbf R} \left((Hp)(x)\right)^2p(x)\,dx =
\frac{\pi^2}{3}\int_{\mathbf R} p(x)^3\,dx, 
\end{equation}  and hence the free Fisher information has an alternative
description: 
$$
\Phi(\mu) = 4 \int_{\mathbf R} \left((Hp)(x)\right)^2 p(x)\,dx = 
4  \int_{\mathbf R} \left((Hp)(x)\right)^2 d\mu(x). 
$$
Here, we should remark that the Hilbert transform is usually defined
with an additional multiple constant $1/\pi$ and
$\int_{\mathbf R} \left((Hp)(x)\right)^2 p(x) dx =
\frac{1}{3} \int_{\mathbf R} p(x)^3 dx$ holds instead of \eqref{integral
formula for Hilbert transf} in this case.
 
Let $Q$ be a real-valued $C^1$ function on ${\mathbf R}$. For each $\mu \in
{\mathcal M}({\mathbf R})$, Biane and Speicher \cite[\S6]{BS} introduced
the {\it relative free Fisher information} $\Phi_Q(\mu)$ of $\mu$
relative to $Q$, and it is defined to be 
\begin{equation}\label{relative free Fisher on R}
\Phi_Q(\mu) := 4 \int_{\mathbf R} \left(\left(Hp\right)(x) - \frac{1}{2}
Q'(x)\right)^2 d\mu(x) 
\end{equation}
when $\mu$ has the density $p = d\mu/dx$ belonging to $L^3({\mathbf R})$;
otherwise to be $+\infty$.

\subsection{Free entropy and free Fisher information for measures on
${\mathbf T}$}
For each $\mu \in {\mathcal M}({\mathbf T})$, the {\it free entropy}
$\Sigma(\mu)$ of $\mu$ is defined in the same manner as in the real line
case; that is, 
$$
\Sigma(\mu) := \iint_{{\mathbf T}^2}
\log\left|\zeta - \eta\right| d\mu(\zeta)\,d\mu(\eta)  
$$
(\cite[\S\S10.7]{V6}, \cite{HP4}). For its justification to be a right
quantity, see \cite[Proposition 10.8]{V6} in relation to the free Fisher
information as well as \cite[Proposition 1.4]{HP4}, \cite{HP2} from the
microstate approach or large deviation principle. 

Assume that $\mu \in {\mathcal M}({\mathbf T})$ has the density $p =
d\mu/d\zeta$ with respect to the Haar probability measure $d\zeta =
d\theta/2\pi$, $\zeta = e^{\sqrt{-1}\theta}$ with $\theta \in [-\pi,\pi)$ and
further that $p$ belongs to the $L^3$-space
$L^3({\mathbf T}) := L^3({\mathbf T},d\zeta)$. As in the real line case, the
Hilbert transform of $p$ 
\begin{equation}\label{Hilbert transf on T}
(Hp)\Bigl(e^{\sqrt{-1}\theta}\Bigr) :=
\lim_{\varepsilon\searrow0} \int_{\varepsilon \leq |t| <\pi}
\frac{p\Bigl(e^{\sqrt{-1}(\theta-t)}\Bigr)}{\tan\bigl(\frac{t}{2}\bigr)}
\,\frac{dt}{2\pi}
\end{equation}
is important. The principle value limit in \eqref{Hilbert transf on T}
exists for a.e.\ (as long as $p\in L^1(\bT)$), and it is known that
$p \in L^q({\mathbf T})$ implies $Hp \in L^q({\mathbf T})$ as well for
each $1<q<\infty$. See \cite[Chapter V]{Ko} for detailed accounts on the
Hilbert transform on ${\mathbf T}$. Following Voiculescu \cite[\S\S8.9]{V6}
we call the quantity
$$
F(\mu) := \int_{\mathbf T} \left((Hp)(\zeta)\right)^2
d\mu(\zeta) = \int_{\mathbf T} \left((Hp)(\zeta)\right)^2
p(\zeta)\,d\zeta
$$
the {\it free Fisher information} of $\mu$. When $\mu$ has no such density
as above, $F(\mu)$ is defined to be $+\infty$. By \cite[Corollary 8.8 and
Definition 8.9]{V6} the free Fisher information can be written as
$$
F(\mu) = \frac{1}{3}\left(-1+\int_{\mathbf T} p(\zeta)^3\,d\zeta\right).
$$

Let $Q$ be a real-valued $C^1$ function on ${\mathbf T}$. As in the case
of measures on $\bR$, for each $\mu \in {\mathcal M}({\mathbf T})$
we define the {\it relative free Fisher information} $F_Q(\mu)$ to be 
\begin{equation}\label{relative free Fisher on T}
F_Q(\mu) := \int_{\mathbf T}\left((Hp)(\zeta) -
Q'(\zeta)\right)^2 d\mu(\zeta) - \left(\int_{\mathbf T} Q'(\zeta)
\,d\mu(\zeta)\right)^2
\end{equation}
when $\mu$ has the density $p = d\mu/d\zeta$ belonging to
$L^3({\mathbf T})$; otherwise to be $+\infty$. Here, $Q'$ means the
derivative of $Q(e^{\i\theta})$ in $\theta$, i.e.,
$Q'(e^{i\theta}) = \frac{d}{d\theta}Q(e^{\i\theta})$. Slight difference
between the two formulas (\ref{relative free Fisher on R}) and (\ref{relative
free Fisher on T}) is worth notice.

\subsection{Large deviations for self-adjoint random matrices}
\label{large deviation theory for selfadj} 
Let $Q$ be a real-valued continuous function on ${\mathbf R}$ such that 
\begin{equation}\label{growth cond}
\lim_{|x|\rightarrow+\infty} |x|\exp(-\varepsilon Q(x)) = 0
\quad \text{for every $\varepsilon > 0$}.
\end{equation}
The {\it weighted energy integral} associated with $Q$ is defined by  
$$
E_Q(\mu) := -\Sigma(\mu) + \int_{\mathbf R} Q(x)\,d\mu(x) \quad
\text{for $\mu \in {\mathcal M}({\mathbf R})$}. 
$$
According to a fundamental result in the theory of weighted potentials
(see \cite[I.1.3]{ST}), there exists a unique $\mu_Q \in {\mathcal
M}({\mathbf R})$ such that 
$$
E_Q(\mu_Q) = \inf\left\{ E_Q(\mu)  :  \mu \in {\mathcal
M}({\mathbf R}) \right\}, 
$$
and $E_Q(\mu_Q)$ is finite (hence so is $\Sigma(\mu_Q)$). Moreover,
$\mu_Q$ is known to be compactly supported. The minimizer $\mu_Q$ is
sometimes called the {\it equilibrium measure} associated with $Q$. Set
$B(Q) := - E_Q\left(\mu_Q\right)$ so that the function 
\begin{equation}\label{good rate function in selfadj case}
-\Sigma(\mu) + \int_{\mathbf R} Q(x)\,d\mu(x) + B(Q)
\quad \text{for $\mu \in {\mathcal M}({\mathbf R})$}
\end{equation}
is non-negative and is zero only when $\mu = \mu_Q$. It is well known that
if $Q(x) = 2x^2/r^2$ with $r>0$, then the equilibrium measure (or the
unique minimizer) $\mu_Q$ is the $(0,r^2/4)$-{\it semicircular distribution}
$\gamma_{0,r}$ (with variance $r^2/4$):
\begin{equation}\label{semicircular}
d\gamma_{0,r}(x) :=
\frac{2}{\pi r^2}\sqrt{r^2 - x^2}\,\chi_{[-r,r]}(x)\,dx. 
\end{equation} 

For each $n \in {\mathbf N}$ define $\lambda_n(Q) \in
{\mathcal M}(M_n^{sa})$, the $n\times n$ {\it self-adjoint random matrix} 
associated with $Q$, by 
$$
d\lambda_n(Q)(A) := \frac{1}{Z_n(Q)} \exp\bigl(-n\Tr_n(Q(A))\bigr)\,dA,
$$
where $dA$ means the ``Lebesgue measure" on $M_n^{sa} \cong {\mathbf
R}^{n^2}$, i.e., 
$$
dA := \prod_{i=1}^n dA_{ii} \prod_{i<j} d\left({\mathrm Re}\,A_{ij}\right)
d\left({\mathrm Im}\,A_{ij}\right) \quad
\text{with $A = \left[A_{ij}\right]$},
$$
$Q(A)$ is the usual functional calculus and $Z_n(Q)$ is a normalization
constant. It is known (see \cite{Me, HP1} for example) that the {\it joint
eigenvalue distribution} on ${\mathbf R}^n$ of $\lambda_n(Q)$ is given as
$$
d\tilde{\lambda}_n(Q)(x_1,\dots,x_n) :=
\frac{1}{\widetilde{Z}_n(Q)}\exp\left(-n\sum_{i=1}^n
Q(x_i)\right) \prod_{i<j} (x_i - x_j)^2 \prod_{i=1}^n dx_i
$$
with a new normalization constant $\widetilde{Z}_n(Q)$. Moreover, the
{\it mean eigenvalue distribution} on ${\mathbf R}$ of $\lambda_n(Q)$ is
defined by
$$
\hat{\lambda}_n(Q) := \int\cdots\int_{{\mathbf R}^n}
\frac{1}{n}\left(\delta_{x_1}+\cdots+\delta_{x_n}\right)
d\tilde{\lambda}_n(Q)(x_1,\dots,x_n). 
$$

In \cite{BG} Ben Arous and Guionnet showed the large deviation principle
for the empirical eigenvalue distribution of the standard self-adjoint
Gaussian random matrix (i.e., $\lambda_n(Q)$ with $Q(x) = x^2/2$). The
following is its slight generalization given in \cite[5.4.3]{HP1}: When
$\left(x_1,\dots,x_n\right)$ is distributed according to
$\tilde{\lambda}_n(Q)$, the {\it empirical eigenvalue distribution} 
\begin{equation}\label{empirical eigenvalue dist}
\frac{1}{n}\left(\delta_{x_1} + \cdots + \delta_{x_n}\right)
\end{equation}
satisfies the large deviation principle in the scale $1/n^2$ and the good
rate function is given by \eqref{good rate function in selfadj case}.
Furthermore, one has $\displaystyle{B(Q) = \lim_{n\rightarrow\infty}
\frac{1}{n^2}\log\widetilde{Z}_n(Q)}$, i.e., 
\begin{equation}\label{B(Q)}
B(Q) = \lim_{n\rightarrow\infty}\frac{1}{n^2}\log\int\cdots\int_{{\mathbf
R}^n}\exp\left(-n\sum_{i=1}^n Q\left(x_i\right)\right)
\prod_{i<j}(x_i - x_j)^2 \prod_{i=1}^n dx_i.
\end{equation}
See \cite{DZ, DS} for general theory of large deviations. Since
$\mu_Q$ is the unique minimizer of \eqref{good rate function in selfadj
case}, the random measure \eqref{empirical eigenvalue dist} converges in the
weak topology to $\mu_Q$ almost surely, and hence 
\begin{equation}\label{weak convergence of mean eigenvalue dist}
\hat{\lambda}_n(Q) \longrightarrow \mu_Q \quad \text{weakly};
\end{equation}
see \cite[p.\ 211]{HP1} and also \cite{De}. From the viewpoint of
the large deviation theory of level-2 (see \cite{DZ, DS}), the function
\eqref{good rate function in selfadj case} can be regarded as a kind of
free analog of the relative entropy with respect to its unique minimizer
$\mu_Q$. Thus, following Biane and Speicher \cite[\S6]{BS} and Biane
\cite[\S3]{Bi}, we call the function \eqref{good rate function in selfadj
case} the {\it relative free entropy} (or {\it modified free entropy}\,) of
$\mu$ relative to $Q$, which is denoted by $\widetilde{\Sigma}_Q(\mu)$;
that is, 
\begin{equation}\label{modified relative free entropy}
\widetilde{\Sigma}_Q(\mu) := -\Sigma(\mu) +
\int_{\mathbf R} Q(x)\,d\mu(x) + B(Q)
\quad\text{for $\mu\in\cM(\bR)$}. 
\end{equation}
We do {\it not} call this the ``free relative entropy" introduced in
\cite{HMP}, a slightly different relative entropy-like quantity
$\Sigma(\mu,\nu)$ for two probability measures in the framework of free
probability. Indeed, the free relative entropy $\Sigma(\mu,\nu)$
for $\mu,\nu\in\cM(\bR)$ is defined as
$$
\Sigma(\mu,\nu):=\iint_{\bR^2}\log|x-y|\,d(\mu-\nu)(x)\,d(\mu-\nu)(y).
$$
But it is known (see \cite[(2.7)]{HMP}) that
$$
\Sigma(\mu,\mu_Q)=\widetilde\Sigma_Q(\mu)
$$
if the support of $\mu$ is included in that of $\mu_Q$.

\subsection{Large deviations for restricted self-adjoint random matrices}
In the course of finding a right free analog of relative entropy, another
random matrix model associated with $Q$ and $R>0$ was introduced in
\cite{HMP}. Here, $Q$ is an arbitrary real-valued continuous function whose
domain includes $[-R,R]$. The self-adjoint random matrix
$\lambda_n(Q;R) \in {\mathcal M}\left(M_n^{sa}\right)$ restricted on a
compact subset $\left\{A \in M_n^{sa} : \Vert A \Vert_{\infty} \leq
R\right\}$ is defined by 
$$
d\lambda_n(Q;R)(A) := \frac{1}{Z_n(Q;R)}\exp\bigl(-n\Tr_n(Q(A))\bigr)
\chi_{\left\{\Vert A \Vert_\infty \leq R\right\}}(A)\,dA 
$$
with a normalization constant $Z_n(Q;R)$. In the above, $\|\cdot\|_\infty$
means the operator norm. The joint eigenvalue distribution supported in
$[-R,R]^n$ of $\lambda_n(Q;R)$ is given as 
$$
\begin{aligned}
&d\tilde{\lambda}_n(Q;R)(x_1,\dots,x_n) \\
&\qquad:= \frac{1}{\widetilde{Z}_n(Q;R)}
\exp\left(-n\sum_{i=1}^n Q(x_i)\right)
\prod_{i<j}(x_i - x_j)^2\prod_{i=1}^n \chi_{[-R,R]}(x_i)\,dx_i
\end{aligned}
$$
with a new normalization constant $\tilde{Z}_n(Q;R)$. Its mean
eigenvalue distribution $\hat{\lambda}_n(Q;R)$ supported in $[-R,R]$ is
defined as in \S\S1.4. As in the case of $\lambda_n(Q)$, the
following large deviation theorem holds: The finite limit
$\displaystyle{B(Q;R) :=
\lim_{n\rightarrow\infty}\frac{1}{n^2}\log\tilde{Z}_n(Q;R)}$ exists, and
when $\left(x_1,\dots,x_n\right)$ is distributed according to
$\tilde{\lambda}_n(Q;R)$, the empirical eigenvalue distribution
\eqref{empirical eigenvalue dist} satisfies the large deviation principle
in the scale $1/n^2$ with the rate function 
\begin{equation}\label{restricted rate function}
-\Sigma(\mu) + \int_{[-R,R]} Q(x)\,d\mu(x) + B(Q;R) \quad 
\text{for $\mu \in {\mathcal M}([-R,R])$}. 
\end{equation}
The proof of this large deviation principle is similar to \cite[5.4.3 and
5.5.1]{HP1} as noticed in \cite{HMP}. In this setting, there also exists a
unique minimizer $\mu_{Q,R} \in {\mathcal M}([-R,R])$ of the rate function
\eqref{restricted rate function}, whose value at $\mu_{Q,R}$ is
zero. If $R > 0$ is chosen so that $\mu_Q$ in \S\S1.4 is
supported in $[-R,R]$, then $\mu_Q = \mu_{Q,R}$ is seen by comparing the two
rate functions, and hence $B(Q) = B(Q;R)$; this assertion is
essentially same as in \cite[Proposition 2.4]{V2} in the single variable case. 

\subsection{Large deviations for special unitary random matrices}
\label{LDP for unitary RM}
Let $Q$ be a real-valued continuous function on $\bT$. Similarly to the real
line case in \S\S1.4, the weighted energy integral
$$
-\Sigma(\mu)+\int_\bT Q(\zeta)\,d\mu(\zeta)
\quad\text{for $\mu\in\cM(\bT)$}
$$
admits a unique minimizer $\mu_Q\in\cM(\bT)$ (or the equilibrium measure
associated with $Q$). Set
$B(Q):=\Sigma(\mu_Q)-\int_\bT Q(\zeta)\,d\mu_Q(\zeta)$. It is known
(\cite{HP2}) that the function
$$
-\Sigma(\mu)+\int_\bT Q(\zeta)\,d\mu(\zeta)+B(Q)
\quad\text{for $\mu\in\cM(\bT)$}
$$
is the rate function of the large deviation for the empirical eigenvalue
distribution of an $n\times n$ {\it unitary random matrix}
$$
d\lambda_n^\U(Q)(U):={1\over Z_n^\U(Q)}
\exp\Bigl(-n\Tr_n(Q(U))\Bigr)\,dU,
$$
where $dU$ is the Haar probability measure on $\U(n)$, $Q(U)$ is
defined via functional calculus and $Z_n^\U(Q)$ is a normalization constant.
Furthermore,
$$
B(Q)=\lim_{n\to\infty}{1\over n^2}\log\int\cdots\int_{\bT^n}
\exp\Biggl(-n\sum_{i=1}^nQ(\zeta_i)\Biggr)
\prod_{1\le i<j\le n}|\zeta_i-\zeta_j|^2\prod_{i=1}^nd\zeta_i
$$
where $d\zeta_i=d\theta_i/2\pi$ for $\zeta_i=e^{\i\theta_i}$. However,
the above unitary random matrix $\lambda_n^\U(Q)$ is not suitable for
our present purpose as will be explained in \S\S1.7. Thus, we need to modify
the above large deviation to the setup of $\SU(n)$.

Now, we begin with the joint eigenvalue distribution of the Haar probability
measure on the special unitary group $\SU(n)$. Note that the $n$ eigenvalues
$\zeta_1,\dots,\zeta_n$ of $U\in\SU(n)$ satisfy $\zeta_1\cdots\zeta_n=1$,
i.e., $\zeta_n=(\zeta_1\cdots\zeta_{n-1})^{-1}$ so that the joint
density must be a permutation-invariant distribution of
$(\zeta_1,\dots,\zeta_{n-1})\in\bT^{n-1}$. The following explicit form
of the density seems a folklore for specialists, and in fact, it is
easily derived from the Weyl integration formula familiar in representation
theory; see \cite[p.\ 104]{Kn} for example.

\begin{lemma}\label{L-1.1}
The joint eigenvalue distribution on $\bT^{n-1}$ of the Haar probability
measure on $\SU(n)$ is
$$
{1\over n!}\prod_{1\le i<j\le n}|\zeta_i-\zeta_j|^2
\prod_{i=1}^{n-1}d\zeta_i
\quad\text{with $\zeta_n=(\zeta_1\cdots\zeta_{n-1})^{-1}$},
$$
or
\begin{eqnarray*}
&&{1\over n!(2\pi)^{n-1}}\prod_{1\le i<j\le n}
\Big|e^{\i\theta_i}-e^{\i\theta_j}\Big|^2\prod_{i=1}^{n-1}d\theta_i \\
&&\hskip3cm
\text{with $\theta_n=-(\theta_1+\cdots+\theta_{n-1})$ {\rm(mod} $2\pi)$}.
\end{eqnarray*}
\end{lemma}

Let $Q$ be a real-valued continuous function on $\bT$. For each $n\in\bN$
define $\lambda_n(Q)\in\cM(\SU(n))$, the $n\times n$ {\it special
unitary random matrix} associated with $Q$, by
\begin{equation}\label{special unitary random matrix}
d\lambda_n^\SU(Q)(U)
:={1\over Z_n^\SU(Q)}\exp\bigl(-n\Tr_n(Q(U))\bigr)\,dU,
\end{equation}
where $dU$ is the Haar probability measure on $\SU(n)$ and
$Z_n^\SU(Q)$ is a normalization constant. By Lemma \ref{L-1.1} the joint
eigenvalue distribution on $\bT^{n-1}$ of $\lambda_n^\SU(Q)$ is given as
\begin{eqnarray*}
&&d\tilde\lambda_n^\SU(Q)(\zeta_1,\dots,\zeta_{n-1})
={1\over\widetilde Z_n^\SU(Q)}
\exp\Biggl(-n\sum_{i=1}^nQ(\zeta_i)\Biggr)
\prod_{1\le i<j\le n}|\zeta_i-\zeta_j|^2\prod_{i=1}^nd\zeta_i \\
&&\hskip8cm{\rm with}
\quad\zeta_n=(\zeta_1\cdots\zeta_{n-1})^{-1}.
\end{eqnarray*}

The next theorem is the large deviation principle for the empirical
eigenvalue distribution of $\lambda_n^\SU(Q)$, whose proof, based on the
explicit form of the density of $\tilde\lambda_n^\SU(Q)$, will be sketched
in Appendix for the convenience of the reader.

\begin{thm}\label{T-A.2}\label{T-1.2}
The finite limit
$B(Q):=\displaystyle\lim_{n\to\infty}{1\over n^2}\log\widetilde Z_n^\SU(Q)$
exists. When $(\zeta_1,\dots,\zeta_{n-1})$ is distributed on $\bT^{n-1}$
according to $\tilde\lambda_n^\SU(Q)$, the empirical distribution  
${1\over n}(\delta_{\zeta_1}+\cdots+\delta_{\zeta_{n-1}}+\delta_{\zeta_n})$
with $\zeta_n=(\zeta_1\cdots\zeta_{n-1})^{-1}$ satisfies the large
deviation principle in the scale $1/n^2$ with the rate function
\begin{equation}\label{unitary rate function}
\widetilde\Sigma_Q(\mu):=-\Sigma(\mu)+\int_\bT Q(\zeta)\,d\mu(\zeta)+B(Q)
\quad\text{for $\mu\in\cM(\bT)$}.
\end{equation}
Furthermore, there exists a unique minimizer $\mu_Q\in\cM(\bT)$ of the rate
function so that $\widetilde\Sigma_Q(\mu_Q)=0$.
\end{thm}

As before, we call the rate function (\ref{unitary rate function}) the
{\it relative free entropy} of $\mu$ with respect to $Q$, which is denoted by
$\widetilde{\Sigma}_Q(\mu)$ as in \eqref{modified relative free entropy}.

\subsection{Ricci curvature tensor of $\SU(n)$}
Let $M$ be a smooth complete Riemannian manifold of dimension $m$, and let
$\Ric(M)$ denote the {\it Ricci curvature tensor} of $M$. For a real-valued
$C^2$ function $\Psi$ on $M$, the {\it Hessian} of $\Psi$ is denoted by
${\mathrm Hess}(\Psi)$. Our arguments in \S3 and \S5 will need to verify
the so-called {\it Bakry and Emery criterion} with a positive constant
$\rho$:
\begin{equation}\label{B&E}
{\mathrm Ric}(M) + {\mathrm Hess}(\Psi) \geq \rho I_m;
\end{equation}
see \cite{BE} and Theorem 2.1 below.

The Ricci curvature tensor of ${\mathrm U}(n)$ is known to be degenerate,
while that of ${\mathrm SU}(n)$ to be of positive constant (see \cite{Mi}, a
nice reference for the topic) and a straightforward computation shows that
the Ricci curvature tensor of $\SU(n)$ with respect to the Riemannian
structure associated with ${\mathrm Tr}_n$ is 
\begin{equation}\label{Ricci of SU(n)} 
{\mathrm Ric}\left({\mathrm SU}(n)\right) = \frac{n}{2}I_{n^2 -1}.
\end{equation}
This is the reason why we have presented Theorem \ref{T-1.2} with use of
$\SU(n)$ instead of $\U(n)$.

\subsection{Differentiability of trace functions} 
A derivative formula as well as the higher order differentiability for a
certain kind of trace functions will be essential in proving the main result
(Theorem \ref{LSI on T}) in \S3. The topic seems rather familiar to
specialists, however we can find no appropriate literature. Here, a lemma is
recorded in a form tailor-made for our later use without full generality.

Let $f(t)$ be a real-valued function on an interval $(a,b)$, and let
$\lambda_1,\lambda_2,\dots$ be distinct points in $(a,b)$. The {\it divided
differences} $f^{[r]}$ for $r=0,1,2,\dots$ are recursively introduced as
follows: $f^{[0]}(\lambda_1):=f(\lambda_1)$ and
$$
f^{[r]}(\lambda_1,\lambda_2,\dots,\lambda_{r+1})
:={f^{[r-1]}(\lambda_1,\lambda_2,\dots,\lambda_r)
-f^{[r-1]}(\lambda_2,\dots,\lambda_r,\lambda_{r+1})
\over\lambda_1-\lambda_{r+1}}.
$$
When $\lambda_i$'s are not necessarily distinct,
$f^{[r]}(\lambda_1,\lambda_2,\dots,\lambda_{r+1})$ can be defined by
continuity as long as $f\in C^r(a,b)$; for example,
$f^{[2]}(\lambda,\lambda)=f'(\lambda)$ and
$f^{[3]}(\lambda,\lambda,\lambda)=f''(\lambda)/2$. See \cite[\S II.2]{Do} for
basic properties of divided differences. Let $A\in M_n^{sa}$, all of whose
eigenvalues are in $(a,b)$, and $A=\sum_{i=1}^l\lambda_iP_i$ be the spectral
decomposition with distinct eigenvalues $\lambda_1,\dots,\lambda_l$ in
$(a,b)$. For each $H_1,H_2,\dots,H_r\in M_n^{sa}$ we define
\begin{eqnarray*}
&&f^{[r]}(A)\circ(H_1,H_2,\dots,H_r) \\
&&\quad:=\sum_{\sigma\in S_r}\sum_{i_1,\dots,i_{r+1}=1}^l
f^{[r]}(\lambda_{i_1},\lambda_{i_2},\dots,\lambda_{i_{r+1}})
P_{i_1}H_{\sigma(1)}P_{i_2}H_{\sigma(2)}\cdots
P_{i_r}H_{\sigma(r)}P_{i_{r+1}},
\end{eqnarray*}
where $S_r$ is the set of all permutations on $\{1,\dots,r\}$. In particular,
note (\cite[V.3.3]{Bh}) that if $f\in C^1(a,b)$ and
$A=U{\rm diag}(\lambda_1,\dots,\lambda_n)U^*$ is a diagonalization, then 
$$
{d\over dt}\bigg|_{t=0}f(A+tH_1)=f^{[1]}(A)\circ H_1
=U\Bigl(\Bigl[f^{[1]}(\lambda_i,\lambda_j)\Bigr]_{ij}\circ U^*H_1U\Bigr)U^*,
$$
where $\circ$ stands for the Schur product. The next lemma can be shown in an
essentially same way as in the proof of \cite[V.3.3]{Bh}.

\begin{lemma}\label{differential lemma}
Let $A,H_1,\dots,H_m\in M_n^{sa}$ and set $G(x):=A+\sum_{k=1}^mx_kH_k$ for
$x=(x_1,\dots,x_m)\allowbreak\in\bR^m$. Let $f$ be a real-valued $C^r$ function
on
$(a,b)$ for some $r\in\bN$. If the eigenvalues of $G(x)$ are in $(a,b)$ for
all $x$ in an open domain $D$ of $\bR^m$, then the function $\Tr_n(f(G(x)))$
is $C^r$ on $D$ and
$$
\begin{aligned}
{\partial^r\over\partial x_{k_1}\partial x_{k_2}\cdots\partial x_{k_r}}
\Tr_n(f(G(x)))
&=\Tr_n\Bigl(f^{[r]}(G(x))\circ(H_{k_1},H_{k_2},\dots,H_{k_r})\Bigr) \\
&=\Tr_n\Bigl(\Bigl((f')^{[r-1]}(G(x))\circ(H_{k_1},\dots,H_{k_{r-1}})
\Bigr)H_{k_r}\Bigr)
\end{aligned}
$$
for all $1\le k_1,k_2,\dots,k_r\le m$ and $x\in D$. In particular,
$$
{\partial\over\partial x_k}\Tr_n(f(G(x)))=\Tr_n(f'(G(x))H_k)
$$
for all $1\le k\le m$ and $x\in D$.
\end{lemma}

\section{Free LSI for measures on $\bR$}
\setcounter{equation}{0}

In this section we will give a supplementary comment to Biane's work
\cite{Bi} on free version of {\it logarithmic Sobolev inequality}
({\it LSI} for short) for measures on ${\mathbf R}$. LSI's were
first interested in constructive quantum field theory, and it was Gross
\cite{G} who first presented in full generality an LSI for Gaussian
measures. Among huge contributions to the topic, Bakry and Emery \cite{BE}
gave a simple ``local" criterion, the so-called Bakry and Emery criterion
(see \eqref{B&E}), for a given measure to satisfy an LSI. Let $M$ be an
$m$-dimensional smooth complete Riemannian manifold with the volume measure
$dx$. The precise statement that Bakry and Emery established is as follows:

\begin{thm}\label{B&E thm} {\rm (Bakry and Emery \cite{BE})} 
Let $\Psi \in C^2\left(M\right)$, and set $d\nu(x) :=
\frac{1}{Z}e^{-\Psi(x)}dx$ with a normalization constant $Z$. Assume that
the Bakry and Emery criterion
${\mathrm Ric}(M) + {\mathrm Hess}(\Psi) \geq \rho I_m$ holds with a
constant $\rho > 0$. Then, for every $\mu\in\cM(M)$ absolutely continuous
with respect to $\nu$ one has
\begin{equation}\label{classical LSI}
S(\mu,\nu) \leq \frac{1}{2\rho}\int_M
\left\Vert\nabla\log\frac{d\mu}{d\nu}\right\Vert^2 d\mu,
\end{equation}
whenever the density $d\mu/d\nu$ is smooth on $M$.
\end{thm}

Recall that the left-hand side of \eqref{classical LSI} is the relative
entropy \eqref{relative entropy}, while the integral in the right-hand side
is nothing but the (classical) {\it relative Fisher information} of $\mu$
relative to $\nu$.

Motivated by and based on this theorem, the following ``free LSI" was shown
by Biane:

\begin{thm}\label{Biane's thm} {\rm (Biane \cite{Bi})}\quad 
Assume that $Q$ is a real-valued $C^1$ function on ${\mathbf R}$ such that
$Q(x)-{\rho\over2}x^2$ is convex on $\bR$ with a constant
$\rho > 0$. Then, for every $\mu\in\cM(\bR)$ one has
\begin{equation}\label{Biane's free LSI on reals} 
\widetilde{\Sigma}_Q(\mu) \leq \frac{1}{2\rho} \Phi_Q(\mu).
\end{equation}
\end{thm}

Obviously, the above convexity assumption of $Q$ is equivalent to
$Q''(x)\ge\rho$ on $\bR$ as long as $Q$ is a $C^2 function$.

When $Q(x) = \rho x^2/2$ with $\rho > 0$, the relative
free entropy $\widetilde{\Sigma}_Q(\mu)$ is given as
$$
\widetilde{\Sigma}_Q(\mu)=
-\Sigma(\mu) + \frac{\rho}{2} \int_{\mathbf R} x^2 d\mu(x)
- \frac{1}{2}\log\rho - \frac{3}{4}
$$
and its minimizer is the $(0,1/\rho)$-semicircular distribution 
$\gamma_{0,2/\sqrt{\rho}}$ (see \eqref{semicircular}). Thus, in this
special case, for any $\mu\in\cM(\bR)$ having the $L^3$-density $p$ and
satisfying $\int_\bR x^2d\mu(x)<+\infty$, the free LSI becomes
\begin{equation}\label{2-3}
-\Sigma(\mu) + \frac{\rho}{2} \int_{\mathbf R} x^2 d\mu(x)
- \frac{1}{2}\log\rho - \frac{3}{4} \leq \frac{1}{2\rho}\biggl(\Phi(\mu)
- 2\rho + \rho^2 \int_{\mathbf R} x^2 d\mu(x)\biggr)
\end{equation} 
because of $2\int_{\mathbf R}((Hp)(x))xp(x)\,dx = 1$. Indeed, notice
$$
\begin{aligned}
&\int_{\mathbf R}\left(\int_{\mathbf R}
{x-t\over(x-t)^2+\eps^2}\,p(t)\,dt\right)xp(x)\,dx \\
&\qquad=\int_{\mathbf R} p(t)\left(\int_{\mathbf R}
\left(1+{t(x-t)-\eps^2\over(x-t)^2+\eps^2}\right)p(x)\,dx\right)dt \\
&\qquad=1-\int_{\mathbf R} tp(t)\left(\int_{\mathbf R}
{t-x\over(t-x)^2+\eps^2}\,p(x)\,dx\right)dt \\
&\phantom{aaaaaaaaaa}-\int_{\mathbf R} p(t)\left(\int_{\mathbf R}
{\eps^2\over(t-x)^2+\eps^2}\,p(x)\,dx\right)dt
\end{aligned}
$$
so that
$$
\begin{aligned}
&2\int_{\mathbf R}\left(\int_{\mathbf R}
{x-t\over(x-t)^2+\eps^2}\,p(t)\,dt\right)xp(x)\,dx \\
&\qquad=1-\int_{\mathbf R} p(t)\left(\int_{\mathbf R}
{\eps^2\over(t-x)^2+\eps^2}\,p(x)\,dx\right)dt\,.
\end{aligned}
$$
Letting $\eps\searrow0$ gives $2\int(Hp(x))xp(x)\,dx=1$ as long as
$p\in L^3(\bR)$ (see \cite[pp.~92--93]{HP1}). The inequality \eqref{2-3} can
be rewritten as
$$
\chi(\mu) \geq -\frac{1}{2\rho}\Phi(\mu) -\frac{1}{2}\log\rho
+ \frac{1}{2}\log2\pi + 1
$$
thanks to the formula \eqref{free entropy chi}. Maximizing the above
right-hand side over
$\rho > 0$ gives Voiculescu's inequality (\cite[Proposition 7.9]{V2}) 
\begin{equation}\label{Voiculescu's LSI}
\chi(\mu) \geq {1\over2} \log\frac{2\pi e}{\Phi(\mu)}. 
\end{equation} 
(The last argument is contained in \cite[\S\S7.2]{BS}.) In this way, the
free LSI in Theorem \ref{Biane's thm} for the functions $Q(x) =
\rho x^2/2$ with $\rho > 0$ is equivalent to the inequality 
\eqref{Voiculescu's LSI}. 

In \cite[Theorem 3.1]{Bi} Biane proved Theorem \ref{Biane's thm} when both
$Q$ and the density of $\mu$ are sufficiently smooth, and the proof of the
extension to the general case was omitted. It may be also worth noting that
Lemma \ref{differential lemma} was implicitly used in \cite{Bi}. The rest of
this section is a supplement to Biane's proof, completing the
proof of Theorem \ref{Biane's thm}

We need the following general technical lemma.

\begin{lemma}\label{limiting B(Q)}
Let $Q$ and $Q_k$, $k\in\bN$, be real-valued continuous functions
on ${\mathbf R}$ satisfying the following two conditions{\rm :}
\begin{itemize}
\item[(a)] $Q_k$ converges to $Q$ uniformly in any finite interval\,{\rm ;} 
\item[(b)] there exists a real-valued continuous function $\widetilde{Q}$ on
${\mathbf R}$ such that
$$
\lim_{|x|\rightarrow +\infty}
|x| \exp\left(-\varepsilon\widetilde{Q}(x)\right)
= 0\quad \text{for every $\varepsilon > 0$} 
$$ 
and $Q_k(x) \geq \widetilde{Q}(x)$ for all $k\in\bN$
{\rm (}so $Q(x) \geq \widetilde{Q}(x)${\rm )}. 
\end{itemize}
Then, the $B\left(Q_k\right)$'s and $B(Q)$ are defined as finite real
numbers {\rm (see \S\S\ref{large deviation theory for selfadj}\,)}, and one
has $\lim_{k\rightarrow\infty} B(Q_k) = B(Q)$. 
\end{lemma}

\proof
By the assumption (b) we can apply the large deviation theorem for
self-adjoint random matrices associated to the given $Q$ and the $Q_k$'s.
Let $\mu_Q$ and $\mu_{Q_k}$ be the equilibrium measures associated with $Q$
and $Q_k$, respectively, and $R > 0$ is chosen so that $\mu_Q$ is supported
in $[-R,R]$. For each $\varepsilon > 0$, thanks to the assumption (a) we
can choose $k_0$ so that $\left|Q_k(x) - Q(x)\right| < \varepsilon$ for all
$x \in [-R,R]$ and for all $k \geq k_0$. Then for $k \geq k_0$ we have
$$
\begin{aligned} 
B(Q) 
&= B\left(Q;R\right) \\
&= \lim_{n\rightarrow\infty}\frac{1}{n^2}
\log\int_{[-R,R]^n}\exp\left(-n\sum_{i=1}^n Q\left(x_i\right)\right)
\prod_{i<j}\left(x_i - x_j\right)^2 \prod_{i=1}^n dx_i \\
&\leq \liminf_{n\rightarrow\infty}\frac{1}{n^2}
\log\int_{[-R,R]^n}\exp\left(-n\sum_{i=1}^n
\left(Q_k\left(x_i\right)+\varepsilon\right)\right)
\prod_{i<j}\left(x_i - x_j\right)^2 \prod_{i=1}^n dx_i \\ 
&\leq  \varepsilon + \liminf_{n\rightarrow\infty}\frac{1}{n^2}
\log\int_{{\mathbf R}^n}\exp\left(-n\sum_{i=1}^n Q_k\left(x_i\right)\right)
\prod_{i<j}\left(x_i - x_j\right)^2 \prod_{i=1}^n dx_i \\
&= \varepsilon + B\left(Q_k\right)
\end{aligned}
$$
so that $B(Q) \leq \liminf_{k\rightarrow\infty}B(Q_k)$ since $\varepsilon$ is
arbitrary. 

In what follows, we will apply some techniques used in \cite[\S5.5]{HP1} and
\cite{HP3}. For $\alpha > 0$ define 
\begin{alignat*}{2} 
F(x,y) &:= -\log|x-y| + \frac{1}{2}\left(Q(x) + Q(y)\right), &\quad   
F_{\alpha}(x,y) &:= \min\left\{F(x,y), \alpha\right\}; \\
F_k(x,y) &:= -\log|x-y| + \frac{1}{2}\left(Q_k(x) + Q_k(y)\right), &\quad   
F_{k,\alpha}(x,y) &:= \min\left\{F_k(x,y), \alpha\right\}.
\end{alignat*}
Note that the double integrals of $F(x,y)$ and $F_k(x,y)$ with respect to
$\mu \in {\mathcal M}({\mathbf R})$ are the weighted energy integrals
$E_Q(\mu)$ and $E_{Q_k}(\mu)$ associated with $Q$ and $Q_k$, respectively.
Since the tightness of $(\mu_{Q_k})$ can be shown as in the proof of
\cite[5.5.3]{HP1}, a subsequence $(\mu_{Q_{k(l)}})$ can be chosen so that
$\mu_{Q_{k(l)}}$ weakly converges to some $\mu_0 \in {\mathcal M}({\mathbf R})$
and 
\begin{eqnarray*}
&&\lim_{l\rightarrow\infty} \iint_{{\mathbf R}^2} F_{k(l)}(x,y)
\,d\mu_{Q_{k(l)}}(x)\,d\mu_{Q_{k(l)}}(y) \\
&&\qquad= \liminf_{k\rightarrow\infty}
\iint_{{\mathbf R}^2} F_k(x,y)\,d\mu_{Q_k}(x)\,d\mu_{Q_k}(y)
=\liminf_{k\to\infty}\,(-B(Q_k)).
\end{eqnarray*}
As in the proof of \cite[5.5.2]{HP1} it is seen that $F_{k,\alpha}(x,y)
\rightarrow F_{\alpha}(x,y)$ uniformly as $k\rightarrow\infty$ for each
$\alpha > 0$. Hence, we have
$$
\begin{aligned}
-B(Q) &\leq \iint_{{\mathbf R}^2} F(x,y)\,d\mu_0(x)\,d\mu_0(y) \\
&= \sup_{\alpha > 0}\iint_{{\mathbf R}^2}
F_{\alpha}(x,y)\,d\mu_0(x)\,d\mu_0(y) \\
&= \sup_{\alpha > 0} \lim_{l\rightarrow\infty}\iint_{{\mathbf R}^2}
F_{k(l),\alpha}(x,y)\,d\mu_{Q_{k(l)}}(x)\,d\mu_{Q_{k(l)}}(y) \\
&\leq \lim_{l\rightarrow\infty}\iint_{{\mathbf R}^2} F_{k(l)}(x,y)
\,d\mu_{Q_{k(l)}}(x)\,d\mu_{Q_{k(l)}}(y) \\
&= \liminf_{k\rightarrow\infty}\left(-B\left(Q_k\right)\right), 
\end{aligned}
$$
where the first inequality comes from that $\mu_Q$ is a minimizer of
$E_Q(\mu)$ with $-B(Q)=E_Q(\mu_Q)$. Thus
$B(Q) \geq \limsup_{k\rightarrow\infty}B(Q_k)$ follows. \qed 

\bigskip 
Now, let us prove Theorem \ref{Biane's thm} for the general case.
Assume that $\mu$ has the density $p = d\mu/dx \in L^3({\mathbf R})$
and moreover that $\Phi_Q(\mu) = 4 \int \left(\left(Hp\right)(x) -
\frac{1}{2}Q'(x)\right)^2 d\mu(x)$ is finite. Since
$Hp \in L^2\left({\mathbf R},\mu\right)$ by the former assumption, the
latter implies $Q' \in L^2\left({\mathbf R},\mu\right)$ as well.

At first, suppose further that $\mu$ is compactly supported. For each
$\varepsilon >0$ choose a non-negative $C^{\infty}$ function
$\phi_{\varepsilon}$ supported in $[-\varepsilon,\varepsilon]$ with
$\int \phi_{\varepsilon}(x)\,dx = 1$, and consider the convolution
$Q_{\varepsilon} := Q * \phi_{\varepsilon}$. Then $Q_\eps$'s are $C^\infty$
functions, and $Q_\eps\to Q$ and $Q_\eps'\to Q'$ uniformly on each finite
interval as $\eps\searrow0$. (The last assertion is seen because
$Q_\eps'=Q'*\phi_\eps$ follows from the $C^1$ of $Q$.) The convexity
assumption of $Q$ means that  
$$
\lambda Q\left(x_1\right) + (1-\lambda) Q\left(x_2\right) - 
Q\left(\lambda x_1 + (1-\lambda)x_2\right) \geq 
\frac{\rho}{2}\lambda(1-\lambda)\left(x_1 - x_2\right)^2
$$
for all $x_1, x_2 \in {\mathbf R}$ and $0 < \lambda < 1$. This implies the
same convexity of $Q_{\varepsilon}$ so that
$Q_{\varepsilon}''(x) \geq \rho$ for all $x \in {\mathbf R}$. 
Define $p_\eps:=p*\phi_\eps$ and $\mu_\eps\in\cM(\bR)$ by
$d\mu_\eps(x):=p_\eps(x)\,dx$.  Moreover, consider
$Q_{\mu_\eps}(x):=2\int_\bR\log|x-y|\,d\mu_\eps(y)$, which is a $C^\infty$
function on $\bR$. Then we have $Q_{\mu_\eps}'(x)=2(Hp_\eps)(x)$ for
a.e.\ $x\in\bR$ (see the proof of Lemma \ref{lemma on Hilbert transf
on T}\,(i) in \S3). Hence, the proof of Theorem \ref{Biane's thm} in
\cite{Bi} implies that
\begin{equation}\label{free LSI for Q_eps}
\widetilde{\Sigma}_{Q_{\varepsilon}}(\mu_\eps) \leq
\frac{1}{2\rho}\Phi_{Q_{\varepsilon}}(\mu_\eps)
\quad\text{for $\eps>0$}.
\end{equation}
Since the convexity assumption of $Q$ implies that
$Q_{\varepsilon}(x) \geq ax^2 + b$ for some $a > 0$ and $b \in {\mathbf R}$,
Lemma \ref{limiting B(Q)} gives
$$
\lim_{\varepsilon\searrow0} B(Q_{\varepsilon}) = B(Q).
$$
Furthermore, notice that $\|p_\eps-p\|_{L^3}\to0$ and hence
$\|Hp_\eps-Hp\|_{L^3}\to0$ as $\eps\searrow0$ so that we get
$$
\begin{aligned} 
\lim_{\varepsilon\searrow0} \int_{\mathbf R}
Q_{\varepsilon}(x)\,d\mu_\eps(x)
&= \int_{\mathbf R} Q(x) d\mu(x),  \\
\lim_{\varepsilon\searrow0} \int_{\mathbf R} \left(\left(Hp_\eps\right)(x)
- \frac{1}{2}Q_{\varepsilon}'(x)\right)^2 d\mu_\eps(x)
&= \int_{\mathbf R} \left(\left(Hp\right)(x)
- \frac{1}{2}Q'(x)\right)^2 d\mu(x).
\end{aligned}
$$
From \eqref{free LSI for Q_eps} and the above convergences together with the
upper semicontinuity of $\Sigma(\mu)$ (see \cite[5.3.2]{HP1}) we have
$$
\widetilde{\Sigma}_Q(\mu) \le \liminf_{\varepsilon\searrow0}
\widetilde{\Sigma}_{Q_{\varepsilon}}(\mu_\eps)
\leq\lim_{\varepsilon\searrow0}\frac{1}{2\rho}
\Phi_{Q_{\varepsilon}}(\mu_\eps)
=\frac{1}{2\rho}\Phi_Q(\mu). 
$$

Next, let us treat the case where $\mu$ is not compactly supported. For
$R>0$ set $d\mu_R(x) :=
\frac{1}{\mu\left([-R,R]\right)}\chi_{[-R,R]}(x)\,d\mu(x)$, whose density
is given by $p_R:=\frac{1}{\mu\left([-R,R]\right)}\chi_{[-R,R]}\,p$. Then,
$\|p_R-p\|_{L^3}\to0$ and $\|Hp_R-Hp\|_{L^3}\to0$ as $R \rightarrow +\infty$
so that 
$$
\lim_{R\rightarrow+\infty}\int_{\mathbf R}
(\left(Hp_R)(x)\right)^2 p_R(x)\,dx
= \int_{\mathbf R} \left((Hp)(x)\right)^2 p(x)\,dx
$$
and
$$
\begin{aligned} 
&\int_{\mathbf R} \left(Q'(x)\right)^2 \left|p_R(x) - p(x)\right| dx \\
&\quad\leq \int_{{\mathbf R}\setminus[-R,R]} \left(Q'(x)\right)^2 p(x)\,dx
+ \left(\frac{1}{\mu\left([-R,R]\right)}-1\right) \int_{\mathbf R}
\left(Q'(x)\right)^2 p(x)\,dx \\
&\quad\longrightarrow 0\quad\text{as $R \to +\infty$}.
\end{aligned}
$$
Furthermore, we have
$$
\begin{aligned} 
&\left|\int_{\mathbf R}(Hp_R)(x)Q'(x)\,d\mu_R(x) -
\int_{\mathbf R}(Hp)(x)Q'(x)\,d\mu(x)\right| \\
&\quad\leq \left\{\int_{\mathbf R} \left((Hp_R)(x)
\left(\frac{p_R(x)}{p(x)} - 1\right)\right)^2\,p(x)\,dx\right\}^{1/2}
\left(\int_{\mathbf R} Q'(x)^2 p(x)\,dx\right)^{1/2} \\
&\phantom{aaaaaa} +
\left(\int_{\mathbf R} \left((Hp_R)(x)
-(Hp)(x)\right)^2 p(x)\,dx\right)^{1/2}
\left(\int_{\mathbf R} Q'(x)^2 p(x)\,dx\right)^{1/2} \\
&\quad\leq \left\{\int_{{\mathbf R}\setminus[-R,R]}
\left((Hp_R)(x)\right)^2 p(x)\,dx
+ \left(\frac{1}{\mu\left([-R,R]\right)} - 1\right)^2
\int_{\mathbf R} \left((Hp_R)(x)\right)^2 p(x)\,dx\right\}^{1/2} \\
&\phantom{aaaaaaaaaaaaaaaaaaaaaaaa}\times
\left(\int_{\mathbf R} Q'(x)^2 p(x)\,dx\right)^{1/2} \\
&\phantom{aaaaa} + \left(\int_{\mathbf R}
\left|(Hp_R)(x)-(Hp)(x)\right|^3 dx\right)^{1/3}
\left(\int_{\mathbf R} p(x)^3\,dx\right)^{1/6}
\left(\int_{\mathbf R} Q'(x)^2 p(x)\,dx\right)^{1/2} \\
&\quad\longrightarrow 0\quad\text{as $R\to +\infty$}.
\end{aligned}
$$
In the above, the first inequality is obtained by the
Cauchy-Schwarz inequality with respect to $d\mu(x) = p(x)dx$ and
the second one is by the H\"{o}lder inequality with respect to $dx$. From
the above convergences we get
\begin{equation}\label{limit Phi_Q}
\lim_{R\rightarrow+\infty}\Phi_Q\left(\mu_Q\right) =
\Phi_Q(\mu).
\end{equation}
On the other hand, we get
\begin{equation}\label{liminf}
\widetilde\Sigma_Q(\mu)\le\liminf_{R\to+\infty}\widetilde\Sigma_Q(\mu_R)
\end{equation}
thanks to the monotone convergence theorem and the upper semicontinuity of
$\Sigma(\mu)$. Therefore, the desired inequality follows from \eqref{limit
Phi_Q}, \eqref{liminf} and the first case of $\mu$ being compactly
supported.\qed   

\section{Free LSI for measures on $\bT$}
\setcounter{equation}{0} 

In this section we will proceed to the free analog of logarithmic Sobolev
inequalities for measures on ${\mathbf T}$. The idea here is essentially same
as Biane's work \cite{Bi} mentioned in \S2. Namely, the free analog arises
as the scaling limit in the scale $1/n^2$ of the classical one
\eqref{classical LSI} on the special unitary group ${\mathrm SU}(n)$. However,
there is an essential difference between his argument and ours; we need full
power of large deviation principle (especially the weak convergence of the
empirical eigenvalue distribution to the equilibrium measure almost surely),
while the weak convergence of the mean eigenvalue distribution is enough in
the proof of \cite[Theorem 3.1]{Bi}.

Let us start with some lemmas.  

\begin{lemma}\label{lemma on derivative} 
Let $Q$ be a harmonic function on a neighborhood of the unit disk\break
$\left\{\zeta \in {\mathbf C} : |\zeta| \leq 1\right\}$. For each $n\in\bN$
and each $U\in\SU(n)$ define $Q(U)$ via the functional calculus and set
$\Psi(U):=\Tr_n(Q(U))$. Then one has
\begin{itemize}
\item[(i)] The function $\Psi(U)$ on $\SU(n)$ is $C^\infty$.
\item[(ii)]
$\nabla\Psi(U)
= \i\left(Q'(U) - \frac{1}{n}{\mathrm Tr}_n(Q'(U))I_n\right)$.
\item[(iii)] If $Q\Bigl(e^{\i t}\Bigr)-{\rho\over2}t^2$ is convex on $\bR$
for some constant $\rho>0$, then $\Hess(\Psi)\ge\rho I_{n^2-1}$.
\end{itemize}
\end{lemma}

\proof
Set $f(t) := Q\Bigl(e^{\sqrt{-1}t}\Bigr)$ for $t\in\bR$, and let
$Y_k := \sqrt{-1}X_k$ with $X_k = X_k^*$, $1 \leq k \leq n^2 -1$, be a
basis of the Lie algebra ${\mathfrak su}(n) = \{ T \in M_n\left({\mathbf
C}\right) : T+T^* = 0,\, {\mathrm Tr}_n(T) = 0\}$
($\cong {\mathbf R}^{n^2 -1}$). For any $U_0 = e^{\sqrt{-1}A_0} \in
{\mathrm SU}(n)$ with $\sqrt{-1}A_0 \in {\mathfrak su}(n)$ and for
$x = \left(x_1,\dots,x_{n^2 -1}\right) \in {\mathbf R}^{n^2 -1}$, we write
$$
\Psi\Biggl(\exp\Biggl(\sqrt{-1}
A_0 + \sum_{k=1}^{n^2 -1}x_k Y_k\Biggr)\Biggr)
= \Tr_n\Biggl(f\Biggl(A_0 + \sum_{k=1}^{n^2 -1}x_k X_k\Biggr)\Biggr).
$$
The $C^\infty$ of $f$ on $\bR$ immediately follows from the assumption of
$Q$. In fact, for each $t_0\in\bR$, the function $f(t_0+t)$ has a power
series expansion for $t$ near $0$. Hence, thanks to Lemma
\ref{differential lemma} we have (i) and
$$
\begin{aligned} 
\nabla\Psi(U_0)
&= \sum_{k=1}^{n^2 -1}{\mathrm Tr}_n(f'(A_0)Y_k)Y_k \\
&= \sum_{k=1}^{n^2 -1}{\mathrm Tr}_n\biggl(\biggl(
f'(A_0) - \frac{1}{n}{\mathrm Tr}_n(f'(A_0))I_n\biggr)Y_k\biggr)Y_k \\
&= \sum_{k=1}^{n^2 -1} \biggl\langle\i
\left(f'(A_0) - \frac{1}{n}{\mathrm Tr}_n(f'(A_0))I_n\right),
Y_k\biggr\rangle_{{\mathrm Tr}_n} Y_k \\
&= \i \left( f'(A_0) - \frac{1}{n}{\mathrm Tr}_n(f'(A_0))I_n\right) \\
&= \i \left( Q'(U_0) - \frac{1}{n}{\mathrm Tr}_n(Q'(U_0))I_n\right),
\end{aligned}
$$
implying (ii).

Set $F(t):=Q\Bigl(e^{\i t}\Bigr)-{\rho\over2}t^2$ for $t\in\bR$.
For any $U_0 = e^{\sqrt{-1}A_0} \in {\mathrm SU}(n)$ with
$\sqrt{-1}A_0 \in {\mathfrak su}(n)$ and for
$\left(x_1,\dots,x_{n^2 -1}\right) \in {\mathbf R}^{n^2 -1}$, we have
$$
\begin{aligned}
&\Psi\Biggl(\exp\Biggl(\i A_0+\sum_{k=1}^{n^2-1}x_kY_k\Biggr)\Biggr) \\
&\qquad=\Tr_n\Biggl(F\Biggl(A_0+\sum_{k=1}^{n^2-1}x_kX_k\Biggr)\Biggr)
+{\rho\over2}\Tr_n\Biggl(\Biggl(A_0+\sum_{k=1}^{n^2-1}
x_kX_k\Biggr)^2\Biggr) \\
&\qquad=\Tr_n\Biggl(F\Biggl(A_0+\sum_{k=1}^{n^2-1}x_kX_k\Biggr)\Biggr)
+{\rho\over2}\Tr_n(A_0^2)+\rho\sum_{k=1}^{n^2-1}\Tr_n(A_0X_k)x_k
+{\rho\over2}\sum_{k=1}^{n^2-1}x_k^2.
\end{aligned}
$$
Since $F(t)$ is convex on $\bR$, it is known (\cite[3.1]{OP}) that
$\Tr_n\bigl(F(A_0+\sum_{k=1}^{n^2-1}x_kX_k)\bigr)$ is convex in
$(x_1,\dots,x_{n^2-1})$ so that (iii) follows.
\qed

\begin{lemma}\label{lemma on Hilbert transf on T}
Assume that $\mu \in {\mathcal M}({\mathbf T})$ has a continuous density
$p = d\mu/d\zeta$ and that $Q_{\mu}(\zeta) := 2 \int_{\mathbf T}
\log\left|\zeta - \eta\right|d\mu(\eta)$ is $C^1$ on $\bT$.
Then one has 
\begin{itemize} 
\item[(i)] $Q_{\mu}'(\zeta) = (Hp)(\zeta)$ for a.e.\ $\zeta \in
{\mathbf T}${\rm ;} 
\item[(ii)] $\int_{\mathbf T} \left((Hp)(\zeta)\right)p(\zeta)\,d\zeta = 0$. 
\end{itemize}
\end{lemma} 

\proof
(i) Let $f$ be an arbitrary $C^1$ function on ${\mathbf T}$. Then we have 
$$
\begin{aligned}
&\int_0^{2\pi}
\frac{d}{d\theta}Q_\mu\Bigl(e^{\sqrt{-1}\theta}\Bigr)
f\Bigl(e^{\sqrt{-1}\theta}\Bigr)\frac{d\theta}{2\pi} \\
&\quad= -\int_0^{2\pi} Q_\mu\Bigl(e^{\sqrt{-1}\theta}\Bigr)
\frac{d}{d\theta}f\Bigl(e^{\sqrt{-1}\theta}\Bigr)\frac{d\theta}{2\pi} \\
&\quad= -\lim_{\varepsilon\searrow0}\int_{|\theta -
t|\geq\varepsilon} 2\log\left|e^{\sqrt{-1}\theta}-e^{\sqrt{-1}t}\right| 
\frac{d}{d\theta}f\Bigl(e^{\sqrt{-1}\theta}\Bigr)
p\Bigl(e^{\sqrt{-1}t}\Bigr)\frac{d\theta\times dt}{(2\pi)^2} \\
&\quad= -\lim_{\varepsilon\searrow0}\int_0^{2\pi}\biggl(
\int_{|\theta - t|\geq\varepsilon}\log\left(2(1-\cos(\theta-t))\right) 
\frac{d}{d\theta}f\Bigl(e^{\sqrt{-1}\theta}\Bigr)
\frac{d\theta}{2\pi}\biggr)p\Bigl(e^{\sqrt{-1}t}\Bigr)\frac{dt}{2\pi}, 
\end{aligned} 
$$
where the second equality is due to the fact that
$\log\left|e^{\sqrt{-1}\theta}-e^{\sqrt{-1}t}\right|
\frac{d}{d\theta}f\Bigl(e^{\sqrt{-1}\theta}\Bigr)$ is bounded above.
Integrating by parts we get
$$
\begin{aligned} 
&\int_{|\theta - t|\geq\varepsilon}\log\left(2(1-\cos(\theta-t))\right) 
\frac{d}{d\theta}f\Bigl(e^{\sqrt{-1}\theta}\Bigr)\frac{d\theta}{2\pi} \\
&\quad= -{\log\left(2\left(1-\cos\varepsilon\right)\right)\over2\pi}
\left(f\Bigl(e^{\sqrt{-1}(t+\varepsilon)}\Bigr) 
- f\Bigl(e^{\sqrt{-1}(t-\varepsilon)}\Bigr)\right)
- \int_{|\theta - t|\geq\varepsilon}
\frac{f\Bigl(e^{\sqrt{-1}\theta}\Bigr)}
{\tan\left(\frac{\theta-t}{2}\right)}\,\frac{d\theta}{2\pi}, 
\end{aligned}
$$
and hence 
$$
\begin{aligned}
&\int_0^{2\pi}\frac{d}{d\theta}Q_\mu\Bigl(e^{\sqrt{-1}\theta}\Bigr)
f\Bigl(e^{\sqrt{-1}\theta}\Bigr)\frac{d\theta}{2\pi} \\
&\quad= \lim_{\varepsilon\searrow0}\Biggl\{
\frac{\log\left(2\left(1-\cos\varepsilon\right)\right)}{2\pi}
\int_0^{2\pi}\Bigl(f\Bigl(e^{\sqrt{-1}(t+\varepsilon)}\Bigr) 
- f\Bigl(e^{\sqrt{-1}(t-\varepsilon)}\Bigr)\Bigr)
\,p\Bigl(e^{\sqrt{-1}t}\Bigr)\frac{dt}{2\pi} \\
&\phantom{aaaaaaaaaaa}+ 
\int_0^{2\pi}\left(\int_{|\theta - t|\geq\varepsilon}
\frac{f\Bigl(e^{\sqrt{-1}\theta}\Bigr)}
{\tan\left(\frac{\theta-t}{2}\right)}\frac{d\theta}{2\pi}\right)
p\Bigl(e^{\sqrt{-1}t}\Bigr)\frac{dt}{2\pi}\Biggr\} \\
&\quad= \lim_{\varepsilon\searrow0} \int_0^{2\pi}\left(
\int_{|\theta - t|\geq\varepsilon}
\frac{p\Bigl(e^{\sqrt{-1}t}\Bigr)}
{\tan\left(\frac{\theta-t}{2}\right)}\,\frac{dt}{2\pi}\right)
f\Bigl(e^{\sqrt{-1}\theta}\Bigr)\frac{d\theta}{2\pi} \\
&\quad= \int_0^{2\pi}(Hp)\Bigl(e^{\sqrt{-1}\theta}\Bigr)
f\Bigl(e^{\sqrt{-1}\theta}\Bigr)\frac{d\theta}{2\pi}.  
\end{aligned}
$$
In the above, the second equality comes from
$\left|f\Bigl(e^{\i(t+\eps)}\Bigr)-f\Bigl(e^{\i(t-\eps)}\Bigr)\right|
=O(\eps)$ uniformly for $t\in[0,2\pi)$, and since we have in particular
$p\in L^2\left({\mathbf T}\right)$, the last one does from the
$L^2$-convergence of the involved principle value integral to
$Hp$ (see \cite[12.8.2\,(2)]{Ed}). Hence, the desired
assertion follows since $f$ is arbitrary. 

(ii) is seen by taking the limit as $\varepsilon \searrow 0$ of 
$$
\begin{aligned}
&\int_0^{2\pi}\left(\int_{|t-\theta|\geq\varepsilon}
\frac{p\Bigl(e^{\sqrt{-1}t}\Bigr)}{\tan\left(\frac{\theta-t}{2}\right)}
\,\frac{dt}{2\pi}\right)
p\Bigl(e^{\sqrt{-1}\theta}\Bigr)\frac{d\theta}{2\pi} \\
&\qquad= - \int_0^{2\pi}\left(\int_{|\theta-t|\geq\varepsilon}
\frac{p\Bigl(e^{\sqrt{-1}\theta}\Bigr)}{\tan\left(\frac{t-\theta}{2}\right)}
\,\frac{d\theta}{2\pi}\right)p\Bigl(e^{\sqrt{-1}t}\Bigr)\frac{dt}{2\pi}
\end{aligned}
$$
thanks to the $L^2$-convergence of the principle value integral as mentioned
above. \qed  

\begin{thm}\label{LSI on T}
Let $Q$ be a real-valued $C^1$ function on ${\mathbf T}$ such that
$Q\Bigl(e^{\sqrt{-1}t}\Bigr) - \frac{\rho}{2}t^2$ is convex on ${\mathbf R}$
with a constant $\rho > -1/2$. Then, for every
$\mu \in {\mathcal M}({\mathbf T})$ one has 
\begin{equation}\label{free LSI on T}
\widetilde{\Sigma}_Q(\mu) \leq \frac{1}{1 + 2\rho}F_Q(\mu). 
\end{equation}
\end{thm}

In the special case where $Q \equiv 0$ and $\rho = 0$,
the above \eqref{free LSI on T} becomes
$$
-\Sigma(\mu) \leq F(\mu)
$$
and the equilibrium measure $\mu_Q$ is the uniform distribution $d\zeta$.  

In particular, the theorem implies that $F_Q(\mu)\ge0$; that is,
$$
\int_{\mathbf T}\left((Hp)(\zeta) -
Q'(\zeta)\right)^2 d\mu(\zeta) \ge \left(\int_{\mathbf T} Q'(\zeta)
\,d\mu(\zeta)\right)^2
$$
for every $\mu\in\cM(\bT)$ under the above assumption of $Q$. 
Also, suppose that the equilibrium measure $\mu_Q$ has a continuous density
and its support is $\bT$; then we have
$Q(\zeta)=2\int_\bT\log|\zeta-\eta|\,d\mu_Q(\eta)$ for all
$\zeta\in\bT$ due to \cite[Theorem I.3.1]{ST} so that Lemma
\ref{lemma on Hilbert transf on T} gives $F_Q(\mu_Q)=0$.

\bigskip\noindent
{\it Proof of Theorem \ref{LSI on T}.}\enspace First, let us assume:
\begin{itemize} 
\item[(a)] $Q$ is harmonic on a neighborhood of the unit disk; 
\item[(b)] $\mu$ has a continuous density $p = d\mu/d\zeta$, and
$Q_{\mu}(\zeta) := 2\int_{\mathbf T} \log|\zeta - \eta|\,d\mu(\eta)$ is
harmonic on a neighborhood of the unit disk.
\end{itemize}
For each $n \in {\mathbf N}$ define $n \times n$ special unitary random
matrices $\lambda_n^{{\mathrm SU}}(Q)$ and
$\lambda_n^{{\mathrm SU}}(Q_{\mu})$ as in \eqref{special unitary random
matrix}, i.e.,
$$
\begin{aligned}
d\lambda_n^{\mathrm SU}(Q)(U) &:= \frac{1}{Z_n^{\mathrm SU}(Q)}
\exp(-n{\mathrm Tr}_n(Q(U)))\,dU, \\
d\lambda_n^{\mathrm SU}(Q_\mu)(U) &:= \frac{1}{Z_n^{\mathrm SU}
(Q_{\mu})}\exp(-n{\mathrm Tr}_n(Q_\mu(U)))\,dU.
\end{aligned}
$$
Let $\tilde{\lambda}_n^{{\mathrm SU}}(Q)$ and
$\tilde{\lambda}_n^{{\mathrm SU}}(Q_\mu)$ be their joint eigenvalue
distributions on $\bT^{n-1}$. Also, let
$\hat{\lambda}_n^{{\mathrm SU}}(Q)$ and
$\hat{\lambda}_n^{{\mathrm SU}}(Q_{\mu})$ be their mean eigenvalue
distributions (see \S\S\ref{LDP for unitary RM}). According to Theorem
\ref{T-1.2}, the empirical eigenvalue distribution of
$\lambda_n^{{\mathrm SU}}(Q_{\mu})$ satisfies the large deviation
principle in the scale $1/n^2$ whose rate functions is
$\widetilde{\Sigma}_{Q_{\mu}}(\mu)$. Moreover, note (\cite[Theorem
I.3.1]{ST}) that the equilibrium measure associated with $Q_\mu$ (or the
minimizer of $\widetilde{\Sigma}_{Q_{\mu}}$) is the given $\mu$. This
large deviation principle guarantees the following facts (i) and (ii),
which will be the key ingredients in our arguments below.
\begin{itemize}
\item[(i)] $\hat{\lambda}_n^{\mathrm SU}(Q_{\mu}) \to \mu$
weakly as $n \rightarrow \infty$;
\item[(ii)] the empirical distribution ${1\over n}\left(\zeta_1 + \cdots +
\zeta_n\right)$ weakly converges to $\mu$ almost surely as
$n \rightarrow \infty$ when $(\zeta_1,\dots,\zeta_{n-1})$ is distributed
according to $\tilde{\lambda}_n^{\mathrm SU}(Q_{\mu})$ and
$\zeta_n=(\zeta_1\cdots\zeta_{n-1})^{-1}$.
\end{itemize} 
Set $\Psi_n(U) := n{\mathrm Tr}_n\left(Q(U)\right)$ for $U \in {\mathrm
SU}(n)$. Lemma \ref{lemma on derivative}\,(iii) and
\eqref{Ricci of SU(n)} verify the Bakry and Emery criterion: 
\begin{equation}\label{BE-criterion}
{\mathrm Ric}({\mathrm SU}(n)) + {\mathrm Hess}(\Psi_n)
\geq  \left(\frac{n}{2} + n\rho\right) I_{n^2 -1}. 
\end{equation}
Thus, by Theorem \ref{B&E thm} due to Bakry and Emery we get 
\begin{equation}\label{classical assertion}
S\bigl(\lambda_n^\SU(Q_{\mu}),\lambda_n^\SU(Q)\bigr) \leq
\frac{1}{2\left(\frac{n}{2} + n\rho\right)} \int_{{\mathrm SU}(n)}
\left\Vert \nabla
\log\frac{d\lambda_n^\SU(Q_{\mu})}{d\lambda_n^\SU(Q)}
\right\Vert_{HS}^2 d\lambda_n^\SU(Q_\mu).
\end{equation}
Notice
\begin{equation}\label{RN-derivative}
\frac{d\lambda_n^{\mathrm SU}(Q_{\mu})}{d\lambda_n^{\mathrm SU}(Q)}(U)
=\frac{\widetilde{Z}_n^{\mathrm SU}(Q)}{\widetilde{Z}_n^{\mathrm SU}(Q_\mu)}
\exp\bigl(-n\Tr_n(Q_\mu(U))+n\Tr_n(Q(U))\bigr),
\quad U \in {\mathrm SU}(n),
\end{equation}
where $\widetilde{Z}_n^{\mathrm SU}(Q)$ and
$\widetilde{Z}_n^{\mathrm SU}(Q_\mu)$ are the normalization constants of the
joint eigenvalue distributions (see \S\S\ref{LDP for unitary RM}). Hence, we
have  
$$
\begin{aligned}
&\frac{1}{n^2}S\bigl(\lambda_n^{\mathrm SU}(Q_{\mu}),
\lambda_n^{\mathrm SU}(Q)\bigr) \\
&\quad=\frac{1}{n^2}\int_{{\mathrm SU}(n)}
\log\frac{d\lambda_n^{\mathrm SU}(Q_{\mu})}
{d\lambda_n^{{\mathrm SU}(Q)}(U)}\,d\lambda_n^{\mathrm SU}(Q_{\mu})(U) \\
&\quad=\frac{1}{n^2}\log\widetilde{Z}_n^{\mathrm SU}(Q)
-\frac{1}{n^2}\log\widetilde{Z}_n^{\mathrm SU}(Q_{\mu}) \\
&\quad\phantom{aaa}-\int_{{\mathrm SU}(n)}\frac{1}{n}{\mathrm
Tr}_n(Q_{\mu}(U))
\,d\lambda_n^{\mathrm SU}(Q_{\mu})(U)
+\int_{{\mathrm SU}(n)}\frac{1}{n}{\mathrm Tr}_n(Q(U))
\,d\lambda_n^{\mathrm SU}(Q_{\mu})(U) \\
&\quad=\frac{1}{n^2}\log\widetilde{Z}_n^{\mathrm SU}(Q)
-\frac{1}{n^2}\log\widetilde{Z}_n^{\mathrm SU}(Q_{\mu}) \\
&\quad\phantom{aaa}-\int_{\mathbf T} Q_{\mu}(\zeta)
\,d\hat{\lambda}_n^{\mathrm SU}(Q_{\mu})(\zeta) +\int_{\mathbf T} Q(\zeta)
\,d\hat{\lambda}_n^{\mathrm SU}(Q_{\mu})(\zeta),
\end{aligned}
$$
and therefore, thanks to (b) and (i) above,
\begin{eqnarray}\label{limit of relative entropy}
&&\lim_{n\to\infty}\frac{1}{n^2}S\bigl(\lambda_n^{\mathrm SU}
(Q_{\mu}),\lambda_n^{\mathrm SU}(Q)\bigr) \nonumber\\
&&\qquad=B(Q)-B(Q_\mu)-\int_\bT Q_\mu(\zeta)\,d\mu(\zeta)
+\int_\bT Q(\zeta)\,d\mu(\zeta)=\widetilde\Sigma_Q(\mu),
\end{eqnarray}
where the last equality comes from that $\mu$ is the minimizer with
$\widetilde{\Sigma}_{Q_{\mu}}(\mu) = 0$, i.e., 
$$
\int_{\mathbf T} Q_{\mu}(\zeta)\,d\mu(\zeta) + B(Q_{\mu}) = \Sigma(\mu).
$$
Therefore, the scaling limit in the scale $1/n^2$ of the left-hand side of
\eqref{classical assertion} becomes the relative free entropy
$\widetilde{\Sigma}_Q(\mu)$. We will seek for the scaling limit in the
scale $1/n^2$ of the right-hand side of \eqref{classical assertion}. By
\eqref{RN-derivative} and Lemma \ref{lemma on derivative}\,(ii), we have 
$$
\begin{aligned}
\nabla\log\frac{d\lambda_n^{\mathrm SU}(Q_{\mu})}
{d\lambda_n^{\mathrm SU}(Q)}(U)
&= -n \nabla\bigl({\mathrm Tr}_n(Q_{\mu}(U))-{\mathrm Tr}_n(Q(U))\bigr) \\
&= -\i\Bigl\{n\bigl(Q_{\mu}'(U) - Q'(U)\bigr) -
\left({\mathrm Tr}_n\bigl(Q_{\mu}'(U)-Q'(U)\bigr)\right)I_n\Bigr\}
\end{aligned}
$$
so that 
$$
\begin{aligned}
&\left\Vert \nabla\log\frac{d\lambda_n^{\mathrm SU}(Q_{\mu})}
{d\lambda_n^{\mathrm SU}(Q)}(U) \right\Vert_{HS}^2 \\
&\qquad = n^2{\mathrm Tr}_n\Bigl(\bigl(Q_{\mu}'(U) - Q'(U)\bigr)^2\Bigr)
- n\Bigl({\mathrm Tr}_n\bigl(Q_{\mu}'(U)-Q'(U)\bigr)\Bigr)^2. 
\end{aligned}
$$
Thus, we get 
$$
\begin{aligned}
&\frac{1}{n^2}\cdot\frac{1}{2\left(\frac{n}{2} + n\rho\right)}
\int_{{\mathrm SU}(n)} \left\Vert \nabla\log
\frac{d\lambda_n^{\mathrm SU}(Q_{\mu})}{d\lambda_n^{\mathrm SU}(Q)}(U)
\right\Vert_{HS}^2 d\lambda_n^{\mathrm SU}(Q_\mu)(U) \\
&\qquad= \frac{1}{1 + 2\rho} \Biggl\{
\int_{{\mathrm SU}(n)} \frac{1}{n}
{\mathrm Tr}_n\Bigl(\bigl(Q_{\mu}'(U) - Q'(U)\bigr)^2\Bigr)
\,d\lambda_n^\SU(Q_\mu)(U) \\
&\phantom{aaaaaaaaaaaaaaaa} -
\int_{{\mathrm SU}(n)} \frac{1}{n^2}
\Bigl({\mathrm Tr}_n\bigl(Q_{\mu}'(U)-Q'(U)\bigr)\Bigr)^2
d\lambda_n^\SU(Q_\mu)(U)\Biggr\}.
\end{aligned}
$$
The above-mentioned fact (i) implies that 
$$
\begin{aligned}
&\int_{{\mathrm SU}(n)} \frac{1}{n}
{\mathrm Tr}_n\Bigl(\bigl(Q_{\mu}'(U) - Q'(U)\bigr)^2\Bigr)
\,d\lambda_n^\SU(Q_\mu)(U) \\
&\qquad= \int_{\mathbf T} \left(Q_{\mu}'(\zeta) -
Q'(\zeta)\right)^2 d\hat{\lambda}_n^{\mathrm SU}(Q_\mu)(\zeta) \\
&\qquad\longrightarrow 
\int_{\mathbf T} \left(Q_{\mu}'(\zeta) - Q'(\zeta)\right)^2 d\mu(\zeta)
\quad\text{as $n \rightarrow \infty$},
\end{aligned}
$$
while the above fact (ii) does that 
$$
\begin{aligned} 
&\int_{{\mathrm SU}(n)} \frac{1}{n^2}\Bigl({\mathrm Tr}_n
\bigl(Q_{\mu}'(U) - Q'(U)\bigr)\Bigr)^2 d\lambda_n^\SU(Q_\mu)(U) \\
&\qquad= \int_{{\mathbf T}^{n-1}} \Biggl(\frac{1}{n}\sum_{i=1}^n
\left(Q_{\mu}'(\zeta_i) - Q'(\zeta_i)\right)\Biggr)^2
d\tilde{\lambda}_n^{\mathrm SU}(Q_{\mu})(\zeta_1,\dots,\zeta_{n-1}) \\
&\hskip7.5cm \text{with $\zeta_n := (\zeta_1\cdots\zeta_{n-1})^{-1}$}
\\ &\qquad\longrightarrow
\left(\int_{\mathbf T} \left(Q_{\mu}'(\zeta) - Q'(\zeta)\right)
d\mu(\zeta)\right)^2\quad\text{as $n \rightarrow \infty$}
\end{aligned}
$$
Thanks to the assumption (b), Lemma \ref{lemma on Hilbert transf on T}
implies that
$$
\begin{aligned}
\left(\int_{\mathbf T} \left(Q_{\mu}'(\zeta) - Q'(\zeta)\right)
d\mu(\zeta)\right)^2 &= \left(\int_{\mathbf T} ((Hp)(\zeta)) p(\zeta)
\,d\zeta - \int_{\mathbf T} Q'(\zeta)\,d\mu(\zeta)\right)^2 \\
&= \left(\int_{\mathbf T} Q'(\zeta)\,d\mu(\zeta)\right)^2
\end{aligned}
$$
so that we get
\begin{equation}\label{limit of Fisher}
\lim_{n\rightarrow\infty} \frac{1}{n^2}\cdot\frac{1}{2\left(\frac{n}{2} +
n\rho\right)} \int_{{\mathrm SU}(n)}
\left\Vert \nabla\log\frac{d\lambda_n^{\mathrm SU}(Q_{\mu})}
{d\lambda_n^{\mathrm SU}(Q)}(U)\right\Vert_{HS}^2 d\lambda_n^\SU(Q_\mu)(U)
= \frac{1}{1 + 2\rho}F_Q(\mu).
\end{equation}
By \eqref{classical assertion}, \eqref{limit of relative entropy} and
\eqref{limit of Fisher} we have shown the desired inequality
\eqref{free LSI on T} under the assumptions (a) and (b). 

Next, let us deal with a general $Q$ as stated in the theorem. Let
$\mu \in {\mathcal M}({\mathbf T})$ with a density
$p = d\mu/d\zeta \in L^3({\mathbf T})$. For each $0 < r < 1$, we consider
the Poisson integrals $Q_r$ and $p_r$ of $Q$ and $p$, respectively; that is,
$$
\begin{aligned}
Q_r\Bigl(e^{\sqrt{-1}\theta}\Bigr) &:= \frac{1}{2\pi}
\int_0^{2\pi} P_r(\theta - t) Q\Bigl(e^{\sqrt{-1}t}\Bigr)\,dt, \\
p_r\Bigl(e^{\sqrt{-1}\theta}\Bigr) &:= \frac{1}{2\pi}
\int_0^{2\pi} P_r(\theta - t) p\Bigl(e^{\sqrt{-1}t}\Bigr)\,dt
\end{aligned}
$$
with the Poisson kernel $P_r(\theta) := (1-r^2)/(1-2r\cos\theta + r^2)$.
Define $\mu_r \in {\mathcal M}({\mathbf T})$ by $d\mu_r(\zeta) :=
p_r(\zeta)d\zeta$.  Then it is plain to see that $Q_r$ satisfies the
assumption (a) and that $\mu_r$ does (b). The convexity assumption of $Q$ in
the theorem means that 
$$
\lambda Q\Bigl(e^{\sqrt{-1}s}\Bigr) +
(1-\lambda)Q\Bigl(e^{\sqrt{-1}t}\Bigr) -
Q\Bigl(e^{\sqrt{-1}(\lambda s + (1-\lambda)t)}\Bigr) \geq
\frac{\rho}{2}\lambda(1-\lambda)(t-s)^2
$$
for all $s,t \in {\mathbf R}$ and $0 < \lambda < 1$. It is easy to check that
each $Q_r$, $0 < r < 1$, satisfies the same convexity assumption so that 
\begin{equation}\label{free LSI on T after regularization} 
\widetilde{\Sigma}_{Q_r}(\mu) \leq \frac{1}{1 + 2\rho} F_{Q_r}(\mu)
\end{equation}
by what we have already shown.
It is known (see \cite{HP2} and also \cite[p.224]{HP1}) that
$\mu_r \rightarrow \mu$ weakly and $\Sigma\left(\mu_r\right) \rightarrow
\Sigma(\mu)$ as $r \nearrow 1$. Moreover, it is known (see \cite[5.3.2]{Ko})
that $\left\Vert Q_r - Q\right\Vert_{\infty} \rightarrow 0$ as
$r \nearrow 1$, where $\Vert \cdot \Vert_{\infty}$ means the uniform norm on
$C({\mathbf T})$. Since it is easily seen that 
$$
\left|\frac{1}{n^2}\log\tilde{Z}_n(Q_r) -
\frac{1}{n^2}\log\tilde{Z}_n(Q)\right| \leq \left\Vert Q_r -
Q\right\Vert_{\infty}, 
$$
we have $B(Q_r) \rightarrow B(Q)$ as $r \nearrow 1$. Therefore, we
get 
$$
\lim_{r\nearrow 1}\widetilde{\Sigma}_{Q_r}(\mu) = \widetilde{\Sigma}_Q(\mu). 
$$
Notice that $\left\Vert p_r - p \right\Vert_{L^3} \rightarrow 0$ and hence
$\left\Vert Hp_r - Hp \right\Vert_{L^3} \rightarrow 0$ as $r \nearrow 1$.
Since $Q$ is a $C^1$ function, $Q_r'$ becomes the Poisson integral of $Q'$
so that $\left\Vert Q_r' - Q' \right\Vert_{\infty} \rightarrow 0$ as
$r \nearrow 1$ as well. These imply that 
$$
\begin{aligned} 
\lim_{r\nearrow 1} F_{Q_r}(\mu) 
&= \lim_{r\nearrow 1} \left\{
\int_{\mathbf T} \left((Hp_r)(\zeta) - Q_r'(\zeta)\right)^2
d\mu_r(\zeta) - \left(\int_{\mathbf T} Q_r'(\zeta)\,d\mu_r(\zeta)\right)^2
\right\} \\
&= \int_{\mathbf T} \left((Hp)(\zeta) - Q'(\zeta)\right)^2 d\mu(\zeta)
-\left(\int_{\mathbf T} Q'(\zeta)d\mu(\zeta)\right)^2 = F_Q(\mu). 
\end{aligned} 
$$
Hence, the desired inequality \eqref{free LSI on T} follows by taking the
limit of \eqref{free LSI on T after regularization}. \qed  

\section{Free TCI for measures on $\bR$}
\setcounter{equation}{0}

The second aim of this paper is to obtain the free analog of
transportation cost inequalities for measures on $\bR$ and on $\bT$. We
deal with probability measures on $\bR$ in this section and those on $\bT$
in the next section. The (classical) transportation cost inequalities
compare the Wasserstein distance with the relative entropy (see
\eqref{relative entropy}) for two probability measures. Let us first recall
the definition of the Wasserstein distance. Let $\cX$ be a Polish space
with a metric $d$. The (quadratic) {\it Wasserstein distance} between
$\mu,\nu\in\cM(\cX)$ is defined by
\begin{equation}\label{F-4.1}
W(\mu,\nu):=\inf_{\pi\in\Pi(\mu,\nu)}
\sqrt{\iint_{\cX\times\cX}{1\over2}d(x,y)^2\,d\pi(x,y)},
\end{equation}
where $\Pi(\mu,\nu)$ denotes the set of all probability measures on
$\cX\times\cX$ with marginals $\mu$ and $\nu$, i.e.,
$\pi(\,\cdot\times\cX)=\mu$ and $\pi(\cX\times\cdot\,)=\nu$. The
Wasserstein distance is sometimes defined with the integral of
$d(x,y)^2$ instead of ${1\over2}d(x,y)^2$. The next lemma is well known
and easy to show.

\begin{lemma}\label{L-4.1}
$W(\mu,\nu)$ is weakly lower semicontinuous in $\mu,\nu\in\cM(\cX)$;
namely, if $\mu_n,\nu_n\in\cM(\cX)$, $\mu_n\to\mu$ and $\nu_n\to\nu$ in
the weak topology, then
$$
W(\mu,\nu)\le\liminf_{n\to\infty}W(\mu_n,\nu_n).
$$
\end{lemma}

In the typical case where $\cX=\bR^n$ and $d(x,y)=\|x-y\|$, the usual
Euclidean metric, let $g_n$ be the standard Gaussian measure, i.e.,
$dg_n(x):=(2\pi)^{-n/2}e^{-\|x\|^2/2}\,dx$ ($dx$ means the Lebesgue
measure on $\bR^n$). The celebrated {\it transportation cost inequality}
({\it TCI} for short) of Talagrand \cite{Ta} is
$$
W(\mu,g_n)\le\sqrt{S(\mu,g_n)},
\qquad\mu\in\cM(\bR^n).
$$
This inequality is a bit extended as follows (see \cite{Le}):

\begin{thm}\label{T-4.2}
Let $\Psi:\bR^n\to\bR$ and assume that $\Psi(x)-{\rho\over2}\|x\|^2$ is
convex on $\bR^n$ with a constant $\rho>0$. If
$d\nu(x):={1\over Z}e^{-\Psi(x)}\,dx\in\cM(\bR^n)$ with a normalization
constant $Z$, then
$$
W(\mu,\nu)\le\sqrt{{1\over\rho}S(\mu,\nu)},
\qquad\mu\in\cM(\bR^n).
$$
\end{thm}

In \cite{OV} Otto and Villani established the interrelation between LSI
and TCI by a technique using partial differential equations. Their result,
combined with Bakry and Emery's LSI (\cite{BE} or Theorem \ref{B&E thm}),
implies the following TCI in a setup on Riemannian manifolds, which will
play a crucial role in deriving our free analog of TCI for measures on
$\bT$. In the theorem, let $M$ be an $m$-dimensional smooth complete
Riemannian manifold equipped with the geodesic distance $d(x,y)$ and the
volume measure $dx$.

\begin{thm}\label{T-4.3}
{\rm (Bakry and Emery \cite{BE} and Otto and Villani \cite{OV})}\quad
Let $\Psi$ be a real-valued $C^2$ function on $M$ and set $d\nu(x):=
{1\over Z}e^{-\Psi(x)}\,dx\in\cM(M)$ with a normalization constant $Z$. If
the Bakry and Emery criterion $\Ric(M)+\Hess(\Psi)\ge\rho I_m$ holds with
a constant $\rho>0$, then
$$
W(\mu,\nu)\le\sqrt{{1\over\rho}S(\mu,\nu)},
\qquad\mu\in\cM(M).
$$
\end{thm}

On the other hand, the following free analog of Talagrand's TCI is shown
by Biane and Voiculescu \cite{BV}. Recall that $\gamma_{0,2}$ is the
standard semicircular measure (see (\ref{semicircular})).

\begin{thm}\label{T-4.4}
{\rm (Biane and Voiculescu \cite{BV})}\quad
For every compactly supported $\mu\in\cM(\bR)$,
\begin{equation}\label{F-4.2}
W(\mu,\gamma_{0,2})\le
\sqrt{-\Sigma(\mu)+\int{x^2\over2}\,d\mu(x)-{3\over4}}.
\end{equation}
\end{thm}

In the rest of this section we will present a new proof of the above free
TCI in a more general situation by using a random matrix technique. In
fact, the classical TCI on the matrix space $M_n^{sa}$ asymptotically
approaches to the free analog when the matrix size goes to $\infty$. The
following is our free TCI for probability measures on $\bR$, where the
relative entropy in the classical TCI is replaced by the relative free
entropy (\ref{modified relative free entropy}).

\begin{thm}\label{T-4.5}
Let $Q$ be a real-valued function on $\bR$. If $Q(x)-{\rho\over2}x^2$ is
convex on $\bR$ with a constant $\rho>0$, then
\begin{equation}\label{F-4.3}
W(\mu,\mu_Q)\le\sqrt{{1\over\rho}
\widetilde\Sigma_Q(\mu)}
\end{equation}
for every compactly supported $\mu\in\cM(\bR)$.
\end{thm}

In particular, when $Q(x)=x^2/2$ and so $\rho=1$, the relative free
entropy $\widetilde\Sigma_Q(\mu)$ is the inside of the square root in
(\ref{F-4.2}) and its minimizer is $\gamma_{0,2}$ so that Theorem
\ref{T-4.5} is a generalization of Theorem \ref{T-4.4}.

The next lemma will play a key role in our proof of the theorem.

\begin{lemma}\label{L-4.5} Let $\tilde\mu,\tilde\nu\in\cM(M_n^{sa})$ and
$\hat\mu,\hat\nu$ be the mean eigenvalue distributions on $\bR$ of
$\tilde\mu,\tilde\nu$, respectively. Then
$$
W(\hat\mu,\hat\nu)\le{1\over\sqrt n}W(\tilde\mu,\tilde\nu),
$$
where $W(\tilde\mu,\tilde\nu)$ is the Wasserstein distance with respect
to the distance induced by the Hilbert-Schmidt norm $\|\cdot\|_{HS}$ on
$M_n^{sa}$.
\end{lemma}

\proof For $A\in M_n^{sa}$ let $\lambda_1(A),\dots,\lambda_n(A)$ be the
eigenvalues of $A$ in increasing order with counting multiplicities. The
mean eigenvalue distribution $\hat\mu$ is written as
$$
\hat\mu=\int_{M_n^{sa}}
{1\over n}\left(\delta_{\lambda_1(A)}+\cdots+\delta_{\lambda_n(A)}\right)
d\tilde\mu(A).
$$
For each $\tilde\pi\in\Pi(\tilde\mu,\tilde\nu)$ define
$\hat\pi\in\cM(\bR\times\bR)$ by
$$
\hat\pi(G):=\iint_{M_n^{sa}\times M_n^{sa}}
{1\over n}\,\#\{i:(\lambda_i(A),\lambda_i(B))\in G\}\,d\tilde\pi(A,B)
$$
for Borel sets $G\subset\bR\times\bR$. Since
$$
\hat\pi(F\times\bR)=\int_{M_n^{sa}}
{1\over n}\,\{\#\{i:\lambda_i(A)\in F\}\,d\tilde\mu(A) =\hat\mu(F)
$$
and similarly $\hat\pi(\bR\times F)=\hat\nu(F)$ for $F\subset\bR$, we get
$\hat\pi\in\Pi(\hat\mu,\hat\nu)$ so that
$$
\begin{aligned}
W(\hat\mu,\hat\nu)^2
&\le\iint_{\bR\times\bR}{1\over2}(x-y)^2\,d\hat\pi(x,y) \\
&=\iint_{M_n^{sa}\times M_n^{sa}}
\Biggl\{\iint_{\bR\times\bR}{1\over2}(x-y)^2
\,d\Biggl({1\over n}\sum_{i=1}^n
\delta_{\lambda_i(A)}\otimes\delta_{\lambda_i(B)}\Biggr)\Biggr\}
\,d\tilde\pi(A,B) \\
&={1\over n}\iint_{M_n^{sa}\times M_n^{sa}}{1\over2}\sum_{i=1}^n
\bigl(\lambda_i(A)-\lambda_i(B)\bigr)^2\,d\tilde\pi(A,B).
\end{aligned}
$$
The famous Lidskii-Wielandt majorization for Hermitian matrices (see
\cite{Bh}) implies that
$$
\sum_{i=1}^n\bigl(\lambda_i(A)-\lambda_i(B)\bigr)^2
\le\sum_{i=1}^n\lambda_i(A-B)^2=\|A-B\|_{HS}^2
$$
for all $A,B\in M_n^{sa}$. Therefore,
$$
W(\hat\mu,\hat\nu)^2\le{1\over n}\iint_{M_n^{sa}\times M_n^{sa}}
{1\over2}\|A-B\|_{HS}^2\,d\tilde\pi(A,B),
$$
and taking the infimum over $\tilde\pi\in\Pi(\tilde\mu,\tilde\nu)$ gives
$W(\hat\mu,\hat\nu)^2\le{1\over n}W(\tilde\mu,\tilde\nu)^2$.\qed

\bigskip
\noindent{\it Proof of Theorem \ref{T-4.5}.}\enspace First, let
$\mu\in\cM(\bR)$ be compactly supported, and suppose that the function
$Q_\mu(x):=2\int\log|x-y|\,d\mu(y)$ is finite and continuous on the whole
$\bR$. Choose $R>0$ so that $\mu$ is supported in $[-R,R]$. For each
$n\in\bN$ consider the $n\times n$ self-adjoint random matrix
$\lambda_n(Q_\mu;R)\in\cM(M_n^{sa})$ supported in
$\{A\in M_n^{sa}:\|A\|_\infty\le R\}$ as well as
$\lambda_n(Q)\in\cM(M_n^{sa})$ (see \S\S1.4 and \S\S1.5). Here,
note that the condition (\ref{growth cond}) is automatically satisfied under
the convexity assumption of $Q$. Since the corresponding large deviation
principle guarantees the weak convergence of the mean eigenvalue
distribution $\hat\lambda_n(Q)$ (resp.\ $\hat\lambda_n(Q_\mu;R)$) to
$\mu_Q$ (resp.\ $\mu$), Lemma \ref{L-4.1} gives
\begin{equation}\label{F-4.4}
W(\mu,\mu_Q)\le\liminf_{n\to\infty}
W\bigl(\hat\lambda_n(Q_\mu;R),\hat\lambda_n(Q)\bigr).
\end{equation}
By Lemma \ref{L-4.5} we get
\begin{equation}\label{F-4.5}
W\bigl(\hat\lambda_n(Q_\mu;R),\hat\lambda_n(Q)\bigr)
\le{1\over\sqrt n}W\bigl(\lambda_n(Q_\mu;R),\lambda_n(Q)\bigr).
\end{equation}
Set $\Psi_n(A):=n\Tr_n(Q(A))$ for $A\in M_n^{sa}$; then
$d\lambda_n(Q)(A)={1\over Z_n(Q)}e^{-\Psi_n(A)}\,dA$. Since
$Q(x)-{\rho\over2}x^2$ is convex on $\bR$, so is
$$
\Psi_n(A)-{\rho n\over2}\|A\|_{HS}^2
=n\Tr_n\Bigl(Q(A)-{\rho\over2}A^2\Bigr)
\quad\text{on $M_n^{sa}$}.
$$
Also, note that $\|\cdot\|_{HS}$ corresponds to the Euclidean norm on
$\bR^{n^2}$ under the isometry
$A=[A_{ij}]\in M_n^{sa}\mapsto\bigl((A_{ii})_{1\le i\le n},\,
(\sqrt2A_{ij})_{i<j}\bigr)\in\bR^{n^2}$. Hence, Theorem \ref{T-4.2} implies
that
\begin{equation}\label{F-4.6}
W\bigl(\lambda_n(Q_\mu;R),\lambda_n(Q)\bigr)
\le\sqrt{{1\over\rho n}S\bigl(\lambda_n(Q_\mu;R),\lambda_n(Q)\bigr)}.
\end{equation}
Similarly to the case of special unitary random matrices in the
proof of Theorem \ref{LSI on T}, since
$$
{d\lambda_n(Q_\mu;R)\over d\lambda_n(Q)}(A)
={\widetilde Z_n(Q)\over\widetilde Z_n(Q_\mu;R)}
\exp\bigl(-n\Tr_n(Q_\mu(A))+n\Tr_n(Q(A))\bigr)
$$
on $(M_n^{sa})_R:=\{A\in M_n^{sa}:\|A\|_\infty\le R\}$, we have
\begin{eqnarray}\label{F-4.7}
&&{1\over n^2}S\bigl(\lambda_n(Q_\mu;R),\lambda_n(Q)\bigr) \nonumber\\
&&\quad ={1\over n^2}\log\widetilde Z_n(Q)-{1\over n^2}
\log\widetilde Z_n(Q_\mu)-\int_{(M_n^{sa})_R}{1\over n}\Tr_n(Q_\mu(A))
\,d\lambda_n(Q_\mu;R)(A) \nonumber\\
&&\quad\qquad\qquad
+\int_{(M_n^{sa})_R}{1\over n}\Tr_n(Q(A))\,d\lambda_n(Q_\mu;R)(A) \nonumber\\
&&\quad \longrightarrow B(Q)-B(\mu;R)-\int_{[-R,R]}
Q_\mu(x)\,d\mu(x)+\int_\bR Q(x)\,d\mu(x)=\widetilde\Sigma_Q(\mu)
\end{eqnarray}
thanks to the fact that $\mu$ is the minimizer of the rate function
(\ref{restricted rate function}) with $Q_\mu$ in place of $Q$. Combining
(\ref{F-4.4})--(\ref{F-4.7}) altogether implies the inequality (\ref{F-4.3})
under the continuity assumption of $Q_\mu(x)$.

Finally, let $\mu\in\cM(\bR)$ be a general compactly supported
measure. By the regularization method in \cite[p.\ 216]{HP1} we can
choose a sequence $\{\mu_k\}$ of measures in $\cM(\bR)$ with compact
supports uniformly bounded such that
$Q_{\mu_k}(x)$ is continuous on $\bR$ for each
$k$, $\mu_k\to\mu$ weakly and $\Sigma(\mu_k)\ge\Sigma(\mu)$ for all
$k$. Hence, by Lemma \ref{L-4.1} and the first case we have
$$
\begin{aligned}
W(\mu,\mu_Q)&\le\liminf_{n\to\infty}W(\mu_k,\mu_Q) \\
&\le\liminf_{k\to\infty}\sqrt{{1\over\rho}\widetilde\Sigma_Q(\mu_k)}
\le\sqrt{{1\over\rho}\widetilde\Sigma_Q(\mu)},
\end{aligned}
$$
completing the proof.\qed

\section{Free TCI for measures on $\bT$}
\setcounter{equation}{0}

In this section we will present the free analog of transportation cost
inequalities for measures on $\bT$. The idea with use of special unitary
random matrices is the same as before. In the following we consider two kinds
of Wasserstein distances between probability measures $\mu,\nu\in\cM(\bT)$.
The one is the Wasserstein distance with respect to the usual metric
$|\zeta-\eta|$,
$\zeta,\eta\in\bT$, and the other is with respect to the geodesic
distance (i.e., the angular distance) on $\bT$. We write
$W_{|\cdot|}(\mu,\nu)$ for the former and $W(\mu,\nu)$ for the
latter. Of course, one has
\begin{equation}\label{F-5.1}
W_{|\cdot|}(\mu,\nu)\le W(\mu,\nu),
\qquad\mu,\nu\in\cM(\bT).
\end{equation}

The next theorem is the free TCI for measures on $\bT$ comparing the
Wasserstein distance with the relative free entropy (\ref{unitary rate
function}).

\begin{thm}\label{T-5.1}
Let $Q$ be a real-valued  function on $\bT$. If there exists a
constant $\rho>-{1\over2}$ such that $Q(e^{\i t})-{\rho\over2}t^2$ is convex
on $\bR$, then
\begin{equation}\label{F-5.2}
W_{|\cdot|}(\mu,\mu_Q)\le W(\mu,\mu_Q)
\le\sqrt{{2\over1+2\rho}\widetilde\Sigma_Q(\mu)}
\end{equation}
for every $\mu\in\cM(\bT)$.
\end{thm}

The special case where $Q\equiv0$ and $\rho=0$ is
$$
W_{|\cdot|}\biggl(\mu,{d\theta\over2\pi}\biggr)\le
W\biggl(\mu,{d\theta\over2\pi}\biggr)\le\sqrt{-2\Sigma(\mu)},
\qquad\mu\in\cM(\bT).
$$

We need the next lemma to prove the theorem. Note that the lemma and
the proof remain valid when $\SU(n)$ is replaced by $\U(n)$.

\begin{lemma}\label{L-5.2}
Let $\tilde\mu,\tilde\nu\in\cM(\SU(n))$ and $W(\tilde\mu,\tilde\nu)$
be the Wasserstein distance between $\tilde\mu,\tilde\nu$ with
respect to the geodesic distance on $\SU(n)$. Let $\hat\mu,\hat\nu$
be the mean eigenvalue distributions on $\bT$ of $\tilde\mu,\tilde\nu$,
respectively. Then
$$
W(\hat\mu,\hat\nu)\le{1\over\sqrt n}W(\tilde\mu,\tilde\nu).
$$
\end{lemma}

\proof
We use the symbol $d$ for the geodesic distance on $\SU(n)$ as well
as for that on $\bT$. Define the optimal matching distance on
$\bT^n$ by
$$
\delta(\zeta,\eta):=\min_{\sigma\in S_n}
\sqrt{\sum_{i=1}^nd(\zeta_i,\eta_{\sigma(i)})^2}
$$
for $\zeta=(\zeta_1,\dots,\zeta_n),
\eta=(\eta_1,\dots,\eta_n)\in\bT^n$.
For $U\in\SU(n)$ let
$\lambda(U):=(\lambda_1(U),\dots,\allowbreak\lambda_n(U))$ denote the
element of $\bT^n$ consisting of the eigenvalues of $U$ with multiplicities
and in counter-clockwise order (i.e.,
$0\le\arg\lambda_1(U)\le\cdots\le\arg\lambda_n(U)<2\pi$). First, we
prove
\begin{equation}\label{F-5.3}
\delta(\lambda(U),\lambda(V))\le d(U,V),
\qquad U,V\in\SU(n).
\end{equation}
For $U,V\in\SU(n)$ let $U(t)$ ($0\le t\le1$) be the geodesic curve in
$\SU(n)$ connecting $U$ and $V$. By dividing the curve into
several small pieces if necessary, we may assume that there is a smooth
curve $A(t)$ ($0\le t\le1$) in $\{A\in M_n^{sa}:\Tr_n(A)=0\}$ such that
$U(t)=e^{\i A(t)}$ for $0\le t\le1$. Let $0=t_0<t_1<\dots<t_K=1$ be
any partition of $A(t)$. For $1\le k\le K$ we have
$$
\begin{aligned}
\delta\bigl(\lambda(U(t_{k-1})),\lambda(U(t_k))\bigr)
&\le\Biggl\{\sum_{i=1}^nd\Bigl(e^{\i\lambda_i(A(t_{k-1}))},
e^{\i\lambda_i(A(t_k))}\Bigr)^2\Biggr\}^{1/2} \\
&\le\Biggl\{\sum_{i=1}^n|\lambda_i(A(t_{k-1}))
-\lambda_i(A(t_k))|^2\Biggr\}^{1/2} \\
&\le\|A(t_{k-1})-A(t_k)\|_{HS} \\
&=d(U(t_{k-1}),U(t_k))+o(t_k-t_{k-1}).
\end{aligned}
$$
In the above, $\lambda_1(A_k),\dots,\lambda_n(A_k)$ are the
eigenvalues of $A_k$ in increasing order, and the third inequality
is due to the Lidskii-Wielandt majorization. Therefore,
$$
\delta(\lambda(U),\lambda(V))\le\sum_{k=1}^K
\delta\bigl(\lambda(U(t_{k-1})),\lambda(U(t_k))\bigr)
\le d(U,V)+o(1)
$$
so that (\ref{F-5.3}) follows because $o(1)\to0$ as
$\max_k(t_k-t_{k-1})\to0$.

Now, for each $U,V\in\SU(n)$ let $\sigma_{U,V}\in S_n$ be such that
$$
\delta(\lambda(U),\lambda(V))
=\Biggl\{\sum_{i=1}^nd\bigl(\lambda_i(U),
\lambda_{\sigma_{U,V}(i)}(V)\bigr)^2\Biggr\}^{1/2}.
$$
Of course, we can let $(U,V)\in\SU(n)\times\SU(n)\mapsto
\sigma_{U,V}\in S_n$ measurable. For every
$\tilde\mu,\tilde\nu\in\cM(\SU(n))$ and
$\tilde\pi\in\Pi(\tilde\mu,\tilde\nu)$, define
$\hat\pi\in\cM(\bT\times\bT)$ by
$$
\hat\pi(G):=\iint_{\SU(n)\times\SU(n)}
{1\over n}\,\#\{i:(\lambda_i(U),\lambda_{\sigma_{U,V}(i)}(V))\in G\}
\,d\tilde\pi(U,V)
$$
for Borel sets $G\subset\bT\times\bT$. Since for $F\subset\bT$
$$
\begin{aligned}
\hat\pi(F\times\bT)&=\int_{\SU(n)}
{1\over n}\,\#\{i:\lambda_i(U)\in F\}\,d\tilde\mu(U)=\hat\mu(F), \\
\hat\pi(\bT\times F)&=\int_{\SU(n)}
{1\over n}\,\#\{i:\lambda_i(V)\in F\}\,d\tilde\nu(V)=\hat\nu(F),
\end{aligned}
$$
we have $\hat\pi\in\Pi(\hat\mu,\hat\nu)$ so that
$$
\begin{aligned}
W(\hat\mu,\hat\nu)^2
&\le\iint_{\bT\times\bT}{1\over2}
d(\zeta,\eta)^2\,d\hat\pi(\zeta,\eta) \\
&={1\over n}\iint_{\SU(n)\times\SU(n)}
{1\over2}\sum_{i=1}^nd\bigl(\lambda_i(U),
\lambda_{\sigma_{U,V}(i)}(V)\bigr)^2\,d\tilde\pi(U,V) \\
&={1\over n}\iint_{\SU(n)\times\SU(n)}{1\over2}
\delta(\lambda(U),\lambda(V))^2\,d\tilde\pi(U,V) \\
&\le{1\over n}\iint_{\SU(n)\times\SU(n)}
{1\over2}d(U,V)^2\,d\tilde\pi(U,V)
\end{aligned}
$$
thanks to (\ref{F-5.3}). This implies
$W(\hat\mu,\hat\nu)^2\le{1\over n}W(\tilde\mu,\tilde\nu)^2$.\qed

\bigskip
\noindent{\it Proof of Theorem \ref{T-5.1}.}\enspace
The first inequality of (\ref{F-5.2}) is obvious as noted in (\ref{F-5.1}).
To prove the second, we first assume:
\begin{itemize}
\item[(a)] $Q$ is harmonic on a neighborhood of the unit disk;
\item[(b)] the function
$Q_\mu(\zeta):=2\int_\bT\log|\zeta-\eta|\,d\mu(\eta)$ is finite and
continuous on $\bT$.
\end{itemize}
For each $n\in\bN$ define $\lambda_n^\SU(Q)$, $\lambda_n^\SU(Q_\mu)$ and
$\hat\lambda_n^\SU(Q)$, $\hat\lambda_n^\SU(Q_\mu)$ as in the proof of
Theorem \ref{LSI on T}. Since $\hat\lambda_n^\SU(Q)\to\mu_Q$ and
$\hat\lambda_n^\SU(Q_\mu)\to\mu$ weakly, Lemma \ref{L-4.1} implies that
\begin{equation}\label{F-5.4}
W(\mu,\mu_Q)\le\liminf_{n\to\infty}
W\bigl(\hat\lambda_n^\SU(Q_\mu),\hat\lambda_n^\SU(Q)\bigr).
\end{equation}
On the other hand, Lemma \ref{L-5.2} gives
\begin{equation}\label{F-5.5}
W\bigl(\hat\lambda_n^\SU(Q_\mu),\hat\lambda_n^\SU(Q)\bigr)\le
{1\over\sqrt n}W\bigl(\lambda_n^\SU(Q_\mu),\lambda_n^\SU(Q)\bigr).
\end{equation}
Furthermore, since the function $\Psi_n(U):=n\Tr_n(Q(U))$ on $\SU(n)$
satisfies the Bakry and Emery criterion \eqref{BE-criterion}, Theorem
\ref{T-4.3} implies that
\begin{equation}\label{F-5.6}
W\bigl(\lambda_n^\SU(Q_\mu),\lambda_n^\SU(Q)\bigr)\le
\sqrt{{2\over n+2n\rho}S\bigl(\lambda_n^\SU(Q_\mu),\lambda_n^\SU(Q)\bigr)}.
\end{equation}
The above (\ref{F-5.4})--(\ref{F-5.6}) and \eqref{limit of relative entropy}
(see also Proposition \ref{P-6.1}\,(1) in \S6) altogether prove the second
inequality of (\ref{F-5.2}) under assumptions (a) and (b).

Next, let $Q$ be as stated in the theorem (hence $Q$ is
continuous on $\bT$) and  $\mu\in\cM(\bT)$ be general. For $0<r<1$ let the
Poisson integrals $Q_r$, $p_r$ and $\mu_r$ be as in the proof of
Theorem \ref{LSI on T}. Since $Q_r$ and $\mu_r$ satisfy (a) and (b) above,
the case already shown implies that
\begin{equation}\label{F-5.7}
W(\mu_r,\mu_{Q_r})\le\sqrt{{2\over1+2\rho}
\widetilde\Sigma_{Q_r}(\mu_r)}.
\end{equation}
Moreover, as in the proof of Theorem \ref{LSI on T}, we have
$\|Q_r-Q\|\to0$, $B(Q_r)\to B(Q)$ and
$\widetilde\Sigma_{Q_r}(\mu_r)\to\widetilde\Sigma_Q(\mu)$ as $r\nearrow1$.
Choose any sequence $0<r(k)<1$ with
$r(k)\to1$ such that $\mu_{Q_{r(k)}}\to\mu_0\in\cM(\bT)$ weakly. By the
upper semicontinuity of $\Sigma(\mu)$, we get
$$
0\le\widetilde\Sigma_Q(\mu_0)
\le\liminf_{k\to\infty}\widetilde\Sigma_{Q_{r(k)}}(\mu_{Q_{r(k)}})=0
$$
so that $\mu_0=\mu_Q$. This shows that $\mu_{Q_r}\to\mu_Q$ weakly as
$r\nearrow1$ and
$$
W(\mu,\mu_Q)\le\liminf_{r\nearrow1}W(\mu_r,\mu_{Q_r})
$$
thanks to Lemma \ref{L-4.1}. Hence, the desired inequality finally follows
by taking the limit of \eqref{F-5.7}.\qed

\section{Concluding remarks}
\setcounter{equation}{0}

In this section we collect some remarks, examples and supplementary results.

\subsection{Use of special orthogonal random matrices}

For a real-valued continuous function $Q$, an $n\times n$ special orthogonal
random matrix $\lambda_n^\SO(Q)$ is defined by
$$
d\lambda_n^\SO(Q)(V)
:={1\over Z_n^\SO(Q)}\exp\Bigl(-{n\over2}\Tr_n(Q(V))\Bigr)\,dV,
$$
where $dV$ is the Haar probability measure on the special orthogonal
group $\SO(n)$. The joint eigenvalue distribution on
$\bT^{n-1}$ of $\lambda_n^\SO(Q)$ is
\begin{eqnarray*}
&&d\tilde\lambda_n^\SO(Q)(\zeta_1,\dots,\zeta_{n-1})
={1\over\widetilde Z_n^\SO(Q)}
\exp\Biggl(-{n\over2}\sum_{i=1}^nQ(\zeta_i)\Biggr)
\prod_{1\le i<j\le n}|\zeta_i-\zeta_j|\prod_{i=1}^nd\zeta_i \\
&&\hskip8cm{\rm with}
\quad\zeta_n=(\zeta_1\cdots\zeta_{n-1})^{-1}.
\end{eqnarray*}
The large deviation is analogous to Theorem \ref{T-1.2}; the rate
function is just ${1\over2}\widetilde\Sigma_Q(\mu)$ and its minimizer is
the same $\mu_Q$. On the other hand, note that the Ricci curvature tensor
of $\SO(n)$ is
$$
\Ric(\SO(n))={n-2\over4}I_{n(n-1)/2},
$$
and the Bakry and Emery criterion in place of \eqref{BE-criterion} is
$$
\Ric(\SO(n))+\Hess(\Psi_n)
\ge\left({n-2\over4}+{n\over2}\rho\right)I_{n(n-1)/2},
$$
where $\Psi_n(V):={n\over2}\Tr_n(Q(V))$ for $V\in\SO(n)$. In this way, a
special orthogonal random matrix model can be used as well to obtain the
free LSI in Theorem \ref{LSI on T} and the free TCI in Theorem \ref{T-5.1}.
Similarly, the free TCI in Theorem \ref{T-4.5} can be shown by using
a real symmetric random matrix model
$$
d\lambda_n^{\rm real}(Q)(T):={1\over Z_n^{\rm real}(Q)}
\exp\Bigl(-{n\over2}\Tr_n(Q(T))\Bigr)\,dT,
$$
where $dT:=\prod_{i\le j}dT_{ij}$ on
$M_n(\bR)^{sa}\cong\bR^{n(n+1)/2}$.

\subsection{Some computations}
Let $Q(x):=\rho x^2/2$ on $\bR$ with $\rho>0$. The equilibrium measure
associated with $Q$ is the semicircular measure $\gamma_{0,2/\sqrt\rho}$.
For $\alpha>0$ we compute
$$
\begin{aligned}
\widetilde\Sigma_Q(\gamma_{0,2/\sqrt\alpha})
&={1\over2}\log\alpha+{\rho\over2\alpha}-{1\over2}\log\rho-{1\over2}, \\
\Phi_Q(\gamma_{0,2/\sqrt\alpha})&={(\alpha-\rho)^2\over\alpha}.
\end{aligned}
$$
Since
$$
\lim_{\alpha\to\rho}{\widetilde\Sigma_Q(\gamma_{0,2/\sqrt\alpha})\over
\Phi_Q(\gamma_{0,2/\sqrt\alpha})}={1\over4\rho},
$$
we notice that the bound $1/2\rho$ in the free LSI \eqref{Biane's thm}
cannot be smaller than $1/4\rho$; however it is unknown whether $1/2\rho$
is the best possible bound or not.

For $2\le\lambda\le\infty$ the equilibrium measure associated with
$Q(\zeta):=-(2/\lambda)\Re\zeta$ on $\bT$ is
\begin{equation}\label{meas on T}
\nu_\lambda:=\left(1+{2\over\lambda}\cos\theta\right){d\theta\over2\pi}
\qquad\text{(with $\nu_\infty={d\theta\over2\pi}$)}
\end{equation}
and $\Sigma(\nu_\lambda)=-1/\lambda^2$ (see \cite[5.3.10]{HP1}). When
$4<\lambda\le\infty$, since $Q\Bigl(e^{\i t}\Bigr)+{1\over\lambda}t^2
={2\over\lambda}\Bigl({t^2\over2}-\cos t\Bigr)$
is convex on $\bR$, the free LSI \eqref{free LSI on T} holds with
$1/(1+2\rho)=\lambda/(\lambda-4)$. For example, for $2\le\alpha\le\infty$ we
compute
$$
\widetilde\Sigma_Q(\nu_\alpha)=\Bigl({1\over\alpha}-{1\over\lambda}\Bigr)^2,
\quad F_Q(\nu_\alpha)=2\Bigl({1\over\alpha}-{1\over\lambda}\Bigr)^2.
$$
Again, the optimality of the bound $1/(1+2\rho)$ in \eqref{free LSI on T} is
unknown.

Concerning the free TCI, it does not seem easy to exactly compute the
Wasserstein distance; in fact, we do not know the exact value of
$W(\gamma_{0,r_1},\gamma_{0,r_2})$ for instance.

\subsection{Classical TCI vs.\ free TCI}
Both classical and free TCI's are formulated in terms of the same (quadratic)
Wasserstein distance for measures, and thus it seems interesting to compare
these two. However, in the case of measures on ${\mathbf R}$, the natural
reference measures are Gaussian (not being compactly supported) in the
classical case, while semicircular (being compactly supported) in the free
case, and hence the question is irrelevant in this case. In the case of the
uniform probability measure $d\theta/2\pi$ on ${\mathbf T}$, our free TCI is
$$
W\biggl(\mu,{d\theta\over2\pi}\biggr) \leq \sqrt{-2\Sigma(\mu)},
\qquad\mu \in {\mathcal M}({\mathbf T}),
$$
while to the authors' best knowledge the sharpest classical TCI is
$$
W\biggl(\mu,{d\theta\over2\pi}\biggr) \leq
\sqrt{S\biggl(\mu,{d\theta\over2\pi}\biggr)},
\qquad\mu \in {\mathcal M}({\mathbf T}).
$$
(The latter inequality is seen as follows. It is known (see \cite[p.94]{Le})
that the ``spectral gap" and ``logarithmic Sobolev constant" are the same
number $1$, and \cite[Theorem 1]{OV} implies the desired inequality.) Now,
if the relative free entropy happens to dominate the (usual) relative
entropy up to a positive constant, then a free TCI would immediately follow
from the classical one. However, this is not, and we indeed have the
following examples: 

(1) For an arbitrary $k \in {\mathbf N}$ and for large $n\in\bN$,
let us choose $k$  disjoint intervals $\left[a_j(n),b_j(n)\right]$,
$1\le j\le k$, in ${\mathbf T} = [0,2\pi)$ whose lengths are all $2\pi/kn$
and whose center points are fixed independently of the choice $n$. Consider
$\mu_k(n) \in {\mathcal M}({\mathbf T})$ whose density is
$\sum_{j=1}^k n\chi_{\left[a_j(n), b_j(n)\right]}$. Then we have 
$$
S\biggl(\mu_k(n),\frac{d\theta}{2\pi}\biggr) = \log n.
$$
On the other hand, by a straightforward computation we see that, for a
sufficiently large $n_0 \in {\mathbf N}$, there are constants $c_k<C_k$
depending only on $k$ such that 
$$
c_k + \frac{\log n}{k} \leq -\Sigma\left(\mu_k(n)\right)
\leq C_k + \frac{\log n}{k} \quad \text{for $n \geq n_0$},
$$ 
and thus 
$$
\frac{-\Sigma\left(\mu_k(n)\right)}
{S\bigl(\mu_k(n),\frac{d\theta}{2\pi}\bigr)}
\to\frac{1}{k}\quad\text{as $n\to\infty$}.
$$
The computation is somewhat similar to a free entropy dimension computation
for single variables; see \cite[Proposition 6.1]{S1} for example. 

(2) For the measure $\nu_{\lambda}$ ($2 < \lambda < \infty$) in
\eqref{meas on T}, with the help of a table on integration formulas, we can
compute
$$
S\biggl(\nu_{\lambda},\frac{d\theta}{2\pi}\biggr) = 
\log\left(\frac{1}{2}\left(1+\sqrt{1-\frac{4}{\lambda^2}}\right)\right)
+ 1 +  \frac{4}{\lambda\sqrt{\lambda^2 - 4}}
- \frac{1}{\sqrt{1 - \frac{4}{\lambda^2}}}, 
$$ 
and hence we get 
$$
\frac{S\bigl(\nu_{\lambda},\frac{d\theta}{2\pi}\bigr)}
{-\Sigma\bigl(\nu_{\lambda}\bigr)}
\longrightarrow 0 \quad \text{as $\lambda \rightarrow \infty$.} 
$$

These examples tell us that the minus free entropy $-\Sigma(\mu)$ cannot be
compared with the relative entropy $S(\mu,d\theta/2\pi)$. 

\subsection{Scaling limit formulas for relative free entropy and
relative free Fisher information}
It seems worthwhile to state some scaling limit formulas given in the proofs
of the main theorems in separate propositions, saying that the relative
entropy and the Fisher information of relevant random matrices
asymptotically converge to the corresponding free analogs for limiting
measures. In fact, the formulas for relative free entropy were essentially
got in \cite{HMP}.

The proof of \eqref{limit of relative entropy} gives (1) of the next
proposition, while that of \eqref{limit of Fisher} does (2) because Lemma
\ref{differential lemma} shows that the derivative formula in Lemma
\ref{lemma on derivative}\,(ii) is still valid for any $U\in\SU(n)$ when $Q$
is a real-valued $C^1$ function on $\bT$. The unitary versions are similar.

\begin{prop}\label{P-6.1}
{\rm (1)} Let $Q$ be a real-valued continuous function on ${\mathbf T}$,
and $\mu \in {\mathcal M}({\mathbf T})$.  If
$Q_{\mu}(\zeta) := 2\int_{\mathbf T}\log|\zeta - \eta|\,d\mu(\eta)$ is
finite and continuous on $\bT$, then
$$
\widetilde{\Sigma}_Q(\mu)
=\lim_{n\rightarrow\infty} \frac{1}{n^2}
S\bigl(\lambda_n^{\mathrm SU}(Q_{\mu}),\lambda_n^{\mathrm SU}(Q)\bigr)
=\lim_{n\rightarrow\infty} \frac{1}{n^2}
S\bigl(\lambda_n^{\mathrm U}(Q_{\mu}),\lambda_n^{\mathrm U}(Q)\bigr).
$$

{\rm (2)}  In addition, if $\mu$ has a continuous density
$d\mu/d\zeta$ and both $Q$ and $Q_\mu$ are $C^1$ functions on $\bT$, then
$$
\begin{aligned}
F_Q(\mu)
&= \lim_{n\rightarrow\infty}\frac{1}{n^3}\int_{{\mathrm SU}(n)}
\left\Vert \nabla\log\frac{d\lambda_n^{\mathrm SU}(Q_{\mu})}
{d\lambda_n^{\mathrm SU}(Q)}(U)\right\Vert_{HS}^2
d\lambda_n^{\mathrm SU}(Q_\mu)(U) \\
&= \lim_{n\rightarrow\infty}\frac{1}{n^3}\int_{{\mathrm U}(n)}
\left\Vert \nabla\log\frac{d\lambda_n^{\mathrm U}\left(Q_{\mu}\right)}
{d\lambda_n^{\mathrm U}(Q)}(U)\right\Vert_{HS}^2
d\lambda_n^{\mathrm U}(Q_\mu)(U). 
\end{aligned}
$$
\end{prop}

Similar limit formulas are given also in the real line case. The formula
in (1) below is \eqref{F-4.7}. The proof of (2) is more or less similar to the
circle case; here the fact that $Q_\mu'(x)=2(Hp)(x)$ for a.e.\ $x\in\bR$ is
needed in place of Lemma \ref{lemma on Hilbert transf on T}\,(i). The
details are left to the reader. Note that the limits in both formulas are
independent of the choice of $R$ such that $\mu$ is supported in $[-R,R]$.
Although the assumption of $Q_\mu$ being $C^1$ on $\bR$ seems rather strong,
yet we have many such examples (see \cite[\S IV.5]{ST}).

\begin{prop}\label{P-6.2}
{\rm (1)} Let $Q$ be a real-valued continuous function on $\bR$ satisfying
\eqref{growth cond}, and $\mu\in\cM(\bR)$ be supported in $[-R,R]$. If
$Q_\mu(x):=2\int_\bR\log|x-y|\,d\mu(x)$ is finite and continuous on $\bR$,
then
$$
\widetilde\Sigma_Q(\mu)=\lim_{n\to\infty}
{1\over n^2}S\bigl(\lambda_n(Q_\mu;R),\lambda_n(Q)\bigr).
$$

{\rm (2)} In addition, if $\mu$ has a continuous density $d\mu/dx$ and both
$Q$ and $Q_\mu$ are $C^1$ functions on $\bR$, then
$$
\Phi_Q(\mu)=\lim_{n\to\infty}{1\over n^3}\int_{(M_n^{sa})_R}
\left\|\nabla\log{d\lambda_n(Q_\mu;R)\over d\lambda_n(Q)}(A)
\right\|_{HS}^2d\lambda_n(Q_\mu;R)(A).
$$
\end{prop}

\subsection{Free LSI for measures on $\bR^+$}
The free LSI \eqref{Biane's thm} is applicable in particular for measures
supported in $\bR^+=[0,\infty)$, but we can also show a different inequality
which might be a proper free LSI in the case where the whole space is
$\bR^+$ instead of $\bR$. Let $\cM_s(\bR)$ be the set of symmetric
probability measures on $\bR$. Consider the bijective transformation
$\mu\in\cM(\bR^+)\mapsto\tilde\mu\in\cM_s(\bR)$ defined as
$$
\mu(F)=\tilde\mu(\{x\in\bR:x^2\in F\})
\quad\text{for $F\subset\bR^+$}.
$$
When $\mu\in\cM(\bR^+)$ has the density $p=d\mu/dx$ on $\bR^+$, the measure
$\tilde\mu$ has the density $\tilde p=d\tilde\mu/dx$ on $\bR$ and
$$
\begin{aligned}
\tilde p(x)&=|x|p(x^2),
\quad\text{$x\in\bR$}; \\
p(x)&=\displaystyle{\tilde p(\sqrt x)\over\sqrt x},
\quad\ \text{$x\in\bR^+$}.
\end{aligned}
$$
 
\begin{lemma}\label{L-6.3}
Let $f$ be a measurable function on $\bR^+$ and set $\tilde f(x):=|x|f(x^2)$
for $x\in\bR$. Then $\tilde f\in L^3(\bR,dx)$ if and only
$f\in L^3(\bR^+,x\,dx)$. If this is the case, then the Hilbert transform
$(Hf)(x)$ exists for a.e.\ $x\in\bR^+$ and $(H\tilde f)(x)=x(Hf)(x^2)$
for a.e.\ $x\in\bR$.
\end{lemma}

\proof
The first assertion is seen because
$\int_\bR|\tilde f(x)|^3\,dx=\int_{\bR^+}x|f(x)|^3\,dx$. Suppose $f\in
L^3(\bR^+,x\,dx)$; then $(H\tilde f)(x)$ exists for a.e.\ $x\in\bR$. For
every $x>0$ and $0<\eps<x^2$ we compute
\begin{eqnarray*}
&&x\Biggl(\int_0^{x^2-\eps}+\int_{x^2+\eps}^\infty\Biggr)
{f(t)\over x^2-t}\,dt \\
&&\quad=\Biggl(\int_0^{x^2-\eps}+\int_{x^2+\eps}^\infty\Biggr)
\biggl({1\over x+\sqrt t}+{1\over x-\sqrt t}\biggr)
{\tilde f(\sqrt t)\over2\sqrt t}\,dt \\
&&\quad=\Biggl(\int_0^{\sqrt{x^2-\eps}}
+\int_{\sqrt{x^2+\eps}}^\infty\Biggr)
\biggl({1\over x+s}+{1\over x-s}\biggr)\tilde f(s)\,ds \\
&&\quad=\Biggl(\int_{-\infty}^{\sqrt{x^2-\eps}}
+\int_{2x-\sqrt{x^2-\eps}}^\infty\Biggr){\tilde f(s)\over x-s}\,ds
-\int_{-\sqrt{x^2+\eps}}^{-\sqrt{x^2-\eps}}
{\tilde f(s)\over x-s}\,ds
+\int_{\sqrt{x^2+\eps}}^{2x-\sqrt{x^2-\eps}}
{\tilde f(s)\over x-s}\,ds.
\end{eqnarray*}
The first of the last three terms is the principal value integral converging
to $(H\tilde f)(x)$ as $\eps\searrow0$ for a.e.\ $x>0$, while the second and
the third terms converge to $0$ as $\eps\searrow0$. Indeed,
$$
\int_{-\sqrt{x^2+\eps}}^{-\sqrt{x^2-\eps}}
\bigg|{\tilde f(s)\over x-s}\bigg|\,ds
\le{1\over x+\sqrt{x^2-\eps}}
\int_{-\sqrt{x^2+\eps}}^{-\sqrt{x^2-\eps}}|\tilde f(s)|\,ds
\longrightarrow0
$$
and
\begin{eqnarray*}
&&\int_{\sqrt{x^2+\eps}}^{2x-\sqrt{x^2-\eps}}
\bigg|{\tilde f(s)\over x-s}\bigg|\,ds \\
&&\qquad
\le\biggl(\int_{-\infty}^\infty\big|\tilde f(s)\big|^3\,ds\biggr)^{1/3}
\Biggl(\int_{\sqrt{x^2+\eps}}^{2x-\sqrt{x^2-\eps}}
{ds\over(s-x)^{3/2}}\Biggr)^{2/3} \\
&&\qquad
=\biggl(\int_{-\infty}^\infty\big|\tilde f(s)\big|^3\,ds\biggr)^{1/3}
{2\over\sqrt\eps}\Bigl(\bigl(\sqrt{x^2+\eps}+x\bigr)^{1/2}
-\bigl(x+\sqrt{x^2-\eps}\bigr)^{1/2}\Bigr) \\
&&\qquad\longrightarrow0\quad\text{as $\eps\searrow0$}.
\end{eqnarray*}
Therefore, we see that $(Hf)(x^2)$ exists and $(H\tilde f)(x)=x(Hf)(x^2)$
for a.e.\ $x>0$. Moreover, we have
$(H\tilde f)(x)=-(H\tilde f)(-x)=x(Hf)(x^2)$ for a.e.\ $x<0$ as well.\qed

\bigskip
Let $Q$ be a real-valued $C^1$ function on $\bR^+$. For each
$\mu\in\cM(\bR^+)$ we define the ``relative free Fisher information"
$\Phi_Q^+(\mu)$ to be
$$
\Phi_Q^+(\mu):=4\int_{\bR^+}x\left((Hp)(x)-{1\over2}Q'(x)\right)^2d\mu(x)
$$
when $\mu$ has the density $p=d\mu/dx$ belonging to $L^3(\bR^+,x\,dx)$;
otherwise to be $+\infty$. In particular, the ``free Fisher information"
$\Phi^+(\mu)$ is defined as $\Phi_Q^+(\mu)$ with $Q\equiv0$, i.e.,
$$
\Phi^+(\mu)=4\int_{\bR^+}x(Hp(x))^2\,d\mu(x).
$$

On the other hand, let $Q$ be a real-valued continuous function on $\bR^+$
such that
$$
\lim_{x\to+\infty}x\exp(-\eps Q(x))=0
\quad\text{for any $\eps>0$}.
$$
We define the ``relative free entropy" $\widetilde\Sigma_Q^+(\mu)$ of
$\mu\in\cM(\bR^+)$ as
$$
\widetilde\Sigma_Q^+(\mu):=-\Sigma(\mu)
+\int_{\bR^+}Q(x)\,d\mu(x)+B^+(Q),
$$
where
$$
B^+(Q):=\lim_{n\to\infty}{1\over n^2}\int\cdots\int_{(\bR^+)^n}
\exp\Biggl(-n\sum_{i=1}^nQ(x_i)\Biggr)
\prod_{i<j}(x_i-x_j)^2\prod_{i=1}^ndx_i.
$$
In fact, similarly to the real line case in \S\S1.4, the function
$\widetilde\Sigma_Q^+(\mu)$ on $\cM(\bR^+)$ is the good rate function of the
large deviation principle for the empirical eigenvalue distribution of the
$n\times n$ positive random matrix
$$
d\lambda_n^+(Q)(A) :=\frac{1}{Z_n^+(Q)}\exp\bigl(-n\Tr_n(Q(A))\bigr)
\chi_{\{A\ge0\}}(A)\,dA. 
$$

\begin{prop}\label{P-6.4}
Let $Q$ be a real-valued convex continuous function on $\bR^+$ such that
$Q$ is $C^1$ on $(0,\infty)$ and $Q'(x)\ge\rho$ for all $x>0$ with a
constant
$\rho>0$. Then, for every $\mu\in\cM(\bR^+)$ one has
\begin{equation}\label{F-6.1}
\widetilde\Sigma_Q^+(\mu)\le{1\over\rho}\Phi_Q^+(\mu).
\end{equation}
\end{prop}

\proof
Define $\widetilde Q(x):={1\over2}Q(x^2)$ for $x\in\bR$; then it is easy to
check that $\widetilde Q$ is a $C^1$-function on $\bR$ and
$Q(x)-{\rho\over2}x^2$ is convex on $\bR$. For each $\mu\in\cM(\bR^+)$ we
can apply Theorem \ref{Biane's thm} to $\tilde\mu\in\cM_s(\bR)$
defined as above so that
$$
\widetilde\Sigma_{\widetilde Q}(\tilde\mu)
\le{1\over2\rho}\Phi_{\widetilde Q}(\tilde\mu).
$$
Now, it suffices to show that
\begin{equation}\label{F-6.3}
\Phi_Q^+(\mu)=\Phi_{\widetilde Q}(\tilde\mu),
\end{equation}
\begin{equation}\label{F-6.4}
\widetilde\Sigma_Q^+(\mu)
=2\widetilde\Sigma_{\widetilde Q}(\tilde\mu).
\end{equation}
To prove \eqref{F-6.3}, we may assume that $\mu$ has the density
$p=d\mu/dx\in L^3(\bR^+,x\,dx)$. Letting
$\tilde p = d\tilde\mu/dx \in L^3(\bR,dx)$, we get by Lemma \ref{L-6.3}
$$
\begin{aligned}
\Phi_Q^+(\mu)
&=4\int_{\bR^+}x\biggl((Hp)(x)-{1\over2}Q'(x)\biggr)^2p(x)\,dx \\
&=8\int_{\bR^+}\biggl((H\tilde p)(x)-{1\over2}
\widetilde Q'(x)\biggr)^2\tilde p(x)\,dx
=\Phi_{\widetilde Q}(\tilde\mu).
\end{aligned}
$$
For every $\mu\in\cM(\bR^+)$ we have $\Sigma(\mu)=2\Sigma(\tilde\mu)$
(see \cite[p.\ 198]{HP1}) and
$\int_{\bR^+}Q(x)\,d\mu(x)=2\int_\bR\widetilde Q(x)\,d\tilde\mu(x)$. For
each $\nu\in\cM(\bR)$, setting $\nu'\in\cM(\bR)$ by $\nu'(F):=\nu(-F)$, we
get $\widetilde\Sigma_{\widetilde Q}((\nu+\nu')/2)
\le\widetilde\Sigma_{\widetilde Q}(\nu)$ by the concavity of free entropy
(see \cite[p.\ 193]{HP1}). These facts show that the equilibrium measure
$\mu_{\widetilde Q}$ associated with $\widetilde Q$ coincides with
$\tilde\mu_Q$ where $\mu_Q$ is the unique minimizer of
$\widetilde\Sigma_Q^+(\mu)$. Therefore, we see that
$B^+(Q)=2B(\widetilde Q)$ and
$\widetilde\Sigma_Q^+(\mu)=2\widetilde\Sigma_{\widetilde Q}(\tilde\mu)$ for
all $\mu\in\cM(\bR^+)$.\qed

\bigskip
In particular, when $Q(x)=\rho x$ on $\bR^+$ with $\rho>0$, note that
$\tilde\mu_Q=\gamma_{0,2/\sqrt\rho}$ for the unique minimizer $\mu_Q$ of
$\Sigma_Q^+(\mu)$, and the inequality \eqref{F-6.1} becomes
$$
-\Sigma(\mu)+\rho\int_{\bR^+}x\,d\mu(x)-\log\rho-{3\over2}
\le{1\over\rho}\biggl(\Phi^+(\mu)-2\rho+\rho^2\int_{\bR^+}x\,d\mu(x)\biggr),
$$
that is,
$$
\chi(\mu)\ge-{1\over\rho}\Phi^+(\mu)-\log\rho+{1\over2}\log2\pi+{5\over4}
$$
as long as $\int_{\bR^+}x\,d\mu(x)<+\infty$.
Maximizing the above right-hand side over $\rho>0$ gives
\begin{equation}\label{F-6.5}
\chi(\mu)\ge{1\over2}\log{2\pi e^{1/2}\over\Phi^+(\mu)^2},
\end{equation}
which also follows from \eqref{Voiculescu's LSI} combined with
$\Sigma(\mu)=2\Sigma(\tilde\mu)$ and $\Phi^+(\mu)=\Phi(\tilde\mu)$.
Notice that $\Phi(\mu)$ in \eqref{Voiculescu's LSI} and $\Phi^+(\mu)^2$ in
\eqref{F-6.5} are not comparable. For example, when $\mu\in\cM(\bR^+)$ has
a density $p(x)=(\alpha+1)x^\alpha\chi_{(0,1]}(x)$ with $\alpha>-1/3$, we
compute
$$
\Phi(\mu)={4(\alpha+1)^3\over3(3\alpha+1)},\quad
\Phi^+(\mu)={4(\alpha+1)^3\over3(3\alpha+2)},
$$
so that $\Phi^+(\mu)^2/\Phi(\mu)$ converges to $0$ as $\alpha\to-1/3$ and
also to $+\infty$ as $\alpha\to+\infty$.

\subsection{Free TCI for measures on $\bR^+$}
Consider the bijective transformation
$\mu\in\cM(\bR^+)\mapsto\hat\mu\in\cM(\bR^+)$ defined as
$$
\mu(F)=\hat\mu(\{x\in\bR^+:x^2\in F\})
\quad\text{for $F\subset\bR^+$}.
$$
The next proposition is a free TCI when the whole space is $\bR^+$.

\begin{prop}
Let $Q$ be a real-valued function on $\bR^+$. If $Q(x^2)-\rho x^2$ is convex
on $\bR$ with a constant $\rho>0$, then
$$
W(\hat\mu,\hat\mu_Q)\le\sqrt{{1\over2\rho}\widetilde\Sigma_Q^+(\mu)}
$$
for every compactly supported $\mu\in\cM(\bR^+)$, where $\mu_Q$ is the
minimizer of $\widetilde\Sigma_Q^+(\mu)$.
\end{prop}

\proof
Let $\widetilde Q(x):={1\over2}Q(x^2)$ for $x\in\bR$ as in the proof of
Proposition \ref{P-6.4}. Since $\tilde\mu_Q$ is the minimizer of
$\widetilde\Sigma_{\widetilde Q}(\nu)$ for $\nu\in\cM(\bR)$, Theorem
\ref{T-4.5} and \eqref{F-6.4} imply that
$$
W(\tilde\mu,\tilde\mu_Q)
\le\sqrt{{1\over\rho}\widetilde\Sigma_{\widetilde Q}(\tilde\mu)}
=\sqrt{{1\over2\rho}\widetilde\Sigma_Q^+(\mu)}.
$$
Hence, it remains to show that
$$
W(\hat\mu,\hat\mu_Q)\le W(\tilde\mu,\tilde\mu_Q).
$$
To prove this, let $\pi\in\Pi(\tilde\mu,\tilde\mu_Q)$ and define
$$
\hat\pi(G):=\pi\bigl(\{(x,y)\in\bR\times\bR:(|x|,|y|)\in G\}\bigr)
$$
for Borel sets $G\subset\bR^+\times\bR^+$. Then we get
$\hat\pi\in\Pi(\hat\mu,\hat\mu_Q)$ so that
$$
\begin{aligned}
W(\hat\mu,\hat\mu_Q)
&\le\int_{\bR^+\times\bR^+}{1\over2}(x-y)^2\,d\hat\pi(x,y)
=\int_{\bR\times\bR}{1\over2}(|x|-|y|)^2\,d\pi(x,y) \\
&\le\int_{\bR\times\bR}{1\over2}(x-y)^2\,d\pi(x,y).
\end{aligned}
$$
This implies the desired inequality.\qed

\bigskip
By replacing $\hat\mu$ by $\mu$, the above free TCI can be rewritten as
\begin{eqnarray*}
&&W(\mu,\hat\mu_Q) \\
&&\quad\le\sqrt{{1\over2\rho}\biggl(
-\iint_{\bR^+\times\bR^+}\log|x^2-y^2|\,d\mu(x)\,d\mu(y)
+\int_{\bR^+}Q(x^2)\,d\mu(x)+B^+(Q)\biggr)}
\end{eqnarray*}
for every compactly supported $\mu\in\cM(\bR^+)$. For example, when
$Q(x)=x$ on $\bR^+$ and $\rho=1$, $\hat\mu_Q$ is the quarter-semicircular
distribution ${1\over\pi}\sqrt{4-x^2}\,\chi_{[0,2]}\,dx$ and $B^+(Q)=-3/2$.

\appendix
\section{Proof of Theorem \ref{T-1.2}}
\renewcommand{\thethm}{\Alph{section}.\arabic{thm}}
\renewcommand{\theequation}{\thesection.\arabic{equation}}
\setcounter{equation}{0}

In the following let us keep the relation
$\zeta_n=(\zeta_1\cdots\zeta_{n-1})^{-1}$. The proof below is essentially
same as that in \cite{HP2}. Set
$$
F(\zeta,\eta):=-\log|\zeta-\eta|+{1\over2}(Q(\zeta)+Q(\eta)).
$$
As in \cite{HP2} it suffices to prove the following inequalities:
\begin{itemize}
\item[(i)]
$$
\limsup_{n\to\infty}{1\over n^2}\log\widetilde Z_n^\SU(Q)
\le-\inf_{\mu\in\cM(\bT)}\iint_{\bT^2}
F(\zeta,\eta)\,d\mu(\zeta)\,d\mu(\eta).
$$
\item[(ii)] For every $\mu\in\cM(\bT)$,
\begin{eqnarray*}
&&\inf_G\biggl[\limsup_{n\to\infty}{1\over n^2}
\log\tilde\lambda_n^\SU(Q)\biggl\{{1\over n}
\bigl(\delta_{\zeta_1}+\cdots+\delta_{\zeta_{n-1}}+\delta_{\zeta_n}\bigr)
\in G\biggr\}\biggr] \\
&&\qquad\le-\iint_{\bT^2}F(\zeta,\eta)\,d\mu(\zeta)\,d\mu(\eta)
-\liminf_{n\to\infty}{1\over n^2}\log\widetilde Z_n^\SU(Q),
\end{eqnarray*}
where $G$ runs over all neighborhoods of $\mu$.
\item[(iii)] For every $\mu\in\cM(\bT)$,
$$
\liminf_{n\to\infty}{1\over n^2}\log\widetilde Z_n^\SU(Q)
\ge-\iint_{\bT^2}F(\zeta,\eta)\,d\mu(\zeta)\,d\mu(\eta).
$$
\item[(iv)] For every $\mu\in\cM(\bT)$,
\begin{eqnarray*}
&&\inf_G\biggl[\liminf_{n\to\infty}{1\over n^2}
\log\tilde\lambda_n^\SU(Q)\biggl\{{1\over n}
\bigl(\delta_{\zeta_1}+\cdots+\delta_{\zeta_{n-1}}+\delta_{\zeta_n}\bigr)
\in G\biggr\}\biggr] \\
&&\qquad\ge-\iint_{\bT^2}F(\zeta,\eta)\,d\mu(\zeta)\,d\mu(\eta)
-\limsup_{n\to\infty}{1\over n^2}\log\widetilde Z_n^\SU(Q),
\end{eqnarray*}
where $G$ is as in (ii).
\end{itemize}

The proofs of the first two are the same as in \cite{HP2}, so we
omit them. To prove (iii) and (iv), we may assume (see \cite{HP2}) that
$\mu$ has a continuous density $f>0$ so that
$\mu=f\Bigl(e^{\i\theta}\Bigr)\,d\theta/2\pi$ and
$\delta\le f(\zeta)\le\delta^{-1}$ on $\bT$ for some $\delta>0$. For
each $n\in\bN$ choose
$$
0=b_0^{(n)}<a_1^{(n)}<b_1^{(n)}<a_2^{(n)}<b_2^{(n)}
<\dots<a_n^{(n)}<b_n^{(n)}=2\pi
$$
such that
$$
{1\over2\pi}\int_0^{a_j^{(n)}}f\Bigl(e^{\i\theta}\Bigr)\,d\theta
={j-{1\over2}\over n},\quad
{1\over2\pi}\int_0^{b_j^{(n)}}f\Bigl(e^{\i\theta}\Bigr)\,d\theta
={j\over n}\,;
$$
hence
\begin{equation}\label{F-A.1}
{\pi\delta\over n}\le b_j^{(n)}-a_j^{(n)}\le{\pi\over n\delta},
\quad
{\pi\delta\over n}\le a_j^{(n)}-b_{j-1}^{(n)}\le{\pi\over n\delta}
\end{equation}
for all $1\le j\le n$. Define
$$
\begin{aligned}
\Delta_n&:=\Bigl\{\Bigl(e^{\i\theta_1},\dots,e^{\i\theta_{n-1}}
\Bigr\}\Bigr):a_j^{(n)}\le\theta_j\le b_j^{(n)},\ 1\le j\le n-1\bigr\}, \\
\Theta_n&:=\bigl\{(\theta_1,\dots,\theta_{n-1}):
a_j^{(n)}\le\theta_j\le b_j^{(n)},\ 1\le j\le n-1\bigr\}, \\
\xi_i^{(n)}&:=\max\Bigl\{Q\Bigl(e^{\i\theta}\Bigr):
a_i^{(n)}\le\theta\le b_i^{(n)}\Bigr\}
\ \,{\rm for}\ \,1\le i\le n-1, \\
d_{ij}^{(n)}&:=\min\Bigl\{\Big|e^{\i s}-e^{\i t}\Big|:
a_i^{(n)}\le s\le b_i^{(n)},\ a_j^{(n)}\le t\le b_j^{(n)}\Bigr\}
\ \,{\rm for}\ \,1\le i,j\le n-1.
\end{aligned}
$$
For every neighborhood $G$ of $\mu$, if $n$ is sufficiently large, then
we have
$$
\Delta_n\subset\biggl\{(\zeta_1,\dots,\zeta_{n-1})\in\bT^{n-1}:
{\delta_{\zeta_1}+\cdots+\delta_{\zeta_n}\over n}\in G\biggr\}
$$
so that with $\theta_n=-(\theta_1+\cdots+\theta_{n-1})$
\begin{eqnarray*}
&&\tilde\lambda_n^\SU(Q)\biggl\{{1\over n}
\bigl(\delta_{\zeta_1}+\cdots+\delta_{\zeta_n}\bigr)\in G\biggr\}
\ge\tilde\lambda_n^\SU(Q)(\Delta_n) \\
&&\qquad={1\over\widetilde Z_n^\SU(Q)(2\pi)^{n-1}}
\int\cdots\int_{\Theta_n}
\exp\Biggl(-n\sum_{i=1}^nQ\Bigl(e^{\i\theta_i}\Bigr)\Biggr) \\
&&\hskip4.5cm\times
\prod_{1\le i<j\le n}\Big|e^{\i\theta_i}-e^{\i\theta_j}\Big|^2
\,d\theta_1\cdots d\theta_{n-1} \\
&&\qquad\ge{1\over\widetilde Z_n^\SU(Q)(2\pi)^{n-1}}
\exp\Biggl(-n\sum_{i=1}^{n-1}\xi_i^{(n)}\Biggr)e^{-nM}
\prod_{1\le i<j\le n-1}(d_{ij}^{(n)})^2 \\
&&\hskip2.5cm\times
\int\cdots\int_{\Theta_n}\prod_{i=1}^{n-1}
\Big|e^{\i\theta_i}-e^{-\i(\theta_1+\cdots+\theta_{n-1})}\Big|^2
\,d\theta_1\cdots d\theta_{n-1},
\end{eqnarray*}
where $M:=\max\{Q(\zeta):\zeta\in\bT\}$. Notice
$$
\Bigl\{\theta_1+\cdots+\theta_{n-1}:
(\theta_1,\dots,\theta_{n-1})\in\Theta_n\Bigr\}
=\bigg[\sum_{i=1}^{n-1}a_i^{(n)},\sum_{i=1}^{n-1}b_i^{(n)}\bigg],
$$
and for $n$ large enough
\begin{equation}\label{F-A.2}
\sum_{i=1}^{n-1}b_i^{(n)}-\sum_{i=1}^{n-1}a_i^{(n)}
\ge{n-1\over n}\pi\delta>{3\pi\over n\delta}.
\end{equation}
From (\ref{F-A.1}) and (\ref{F-A.2}) we can choose an interval
$[\alpha,\beta]\subset\Bigl[\sum_{i=1}^{n-1}a_i^{(n)},
\sum_{i=1}^{n-1}b_i^{(n)}\Bigr]$ such that
$\beta-\alpha=\pi\delta/n^2$ and
$$
[-\beta,-\alpha]\subset
\biggl[b_{k-1}^{(n)}+{\pi\delta\over n^2},
a_k^{(n)}-{\pi\delta\over n^2}\biggr]\ \ ({\rm mod}\ 2\pi)
$$
for some $1\le k\le n$. Then, there exist subintervals
$[\alpha_i,\beta_i]\subset\Bigl[a_i^{(n)},b_i^{(n)}\Bigr]$, $1\le i\le n-1$,
such that
$$
\beta_i-\alpha_i={\pi\delta\over n^2(n-1)},\quad
\sum_{i=1}^{n-1}\alpha_i=\alpha,\quad
\sum_{i=1}^{n-1}\beta_i=\beta,
$$
and hence
\begin{eqnarray*}
&&\int\cdots\int_{\Theta_n}\prod_{i=1}^{n-1}
\Big|e^{\i\theta_i}-e^{-\i(\theta_1+\cdots+\theta_{n-1})}\Big|^2
\,d\theta_1\cdots d\theta_{n-1} \\
&&\qquad\ge\int_{\alpha_1}^{\beta_1}\cdots\int_{\alpha_{n-1}}^{\beta_{n-1}}
\Big|e^{\i\theta_i}-e^{-\i(\theta_1+\cdots+\theta_{n-1})}\Big|^2
\,d\theta_1\cdots d\theta_{n-1} \\
&&\qquad\ge\biggl({2\delta\over n^2}\biggr)^{2(n-1)}
\biggl({\pi\delta\over n^2(n-1)}\biggr)^{n-1}.
\end{eqnarray*}
Therefore, for sufficiently large $n$, we get
\begin{eqnarray*}
&&\tilde\lambda_n^\SU(Q)\biggl\{{1\over n}
\bigl(\delta_{\zeta_1}+\cdots+\delta_{\zeta_n}\bigr)\in G\biggr\} \\
&&\qquad\ge{(2\delta^3)^{n-1}\over\widetilde Z_n^\SU(Q)n^{7(n-1)}}
\exp\Biggl(-n\sum_{i=1}^{n-1}\xi_i^{(n)}\Biggr)
\prod_{1\le i<j\le n-1}\Bigl(d_{ij}^{(n)}\Bigr)^2.
\end{eqnarray*}
Since
\begin{eqnarray*}
&&\lim_{n\to\infty}{2\over n^2}\sum_{1\le i<j\le n-1}\log d_{ij}^{(n)} \\
&&\qquad={1\over(2\pi)^2}\int_0^{2\pi}\int_0^{2\pi}
f\Bigl(e^{\i s}\Bigr)f\Bigl(e^{\i t}\Bigr)
\log\Big|e^{\i s}-e^{\i t}\Big|\,ds\,dt \\
&&\qquad=\iint_{\bT^2}\log|\zeta-\eta|\,d\mu(\zeta)\,d\mu(\eta)
\end{eqnarray*}
as well as
$$
\lim_{n\to\infty}{1\over n}\sum_{i=1}^{n-1}\xi_i^{(n)}
={1\over2\pi}\int_0^{2\pi}Q\Bigl(e^{\i s}\Bigr)f\Bigl(e^{\i s}\Bigr)\,ds
=\int_\bT Q(\zeta)\,d\mu(\zeta),
$$
we have
$$
\begin{aligned}
0&\ge\limsup_{n\to\infty}{1\over n^2}\log\tilde\lambda_n^\SU(Q)
\biggl\{{1\over n}
\bigl(\delta_{\zeta_1}+\cdots+\delta_{\zeta_n}\bigr)\in G\biggr\} \\
&\ge-\iint_{\bT^2}F(\zeta,\eta)\,d\mu(\zeta)\,d\mu(\eta)
-\liminf_{n\to\infty}{1\over n^2}\log\widetilde Z_n^\SU(Q)
\end{aligned}
$$
and
\begin{eqnarray*}
&&\liminf_{n\to\infty}{1\over n^2}\log\tilde\lambda_n^\SU(Q)
\biggl\{{1\over n}
\bigl(\delta_{\zeta_1}+\cdots+\delta_{\zeta_n}\bigr)\in G\biggr\} \\
&&\qquad\ge-\iint_{\bT^2}F(\zeta,\eta)\,d\mu(\zeta)\,d\mu(\eta)
-\limsup_{n\to\infty}{1\over n^2}\log\widetilde Z_n^\SU(Q).
\end{eqnarray*}
These imply (iii) and (iv).\qed

\end{document}